\documentclass[a4paper]{article}
\begin{document}
\title {The binary method of integer decomposition and its application}
\author{Puyun Gao}
\date{ }
\maketitle
\begin{abstract}
This paper presents an integer decomposition method. The method
first writes an integer as a polynomial with 2 as variable that its
coefficients are zero or one. Then, suppose that an integer is
decomposed into product of such two polynomials that their orders
are greater than one. By comparing the coefficients of two side, we
obtain a system of the multivariate quadratic equations. Therefore,
the integer decomposition problem is transformed into the problem of
solving a system of the multivariate quadratic equations.
Furthermore, the method transforms integer decomposition problem into
solving the system of linear equations which the number of variables will
not exceed the logarithm with the base of 2 of the decomposed number.
The advantage of this method is direct without any guess which is not available in traditional methods.
\end{abstract}
\section{Introduction}
Integer decomposition is closely related to many famous problems,
such as P and NP problem, and RSA cryptosystem cryptography design
and decryption([1-2]).
\par To break down the RSA cryptosystem is to decompose a large
integer $M$ into the product of two prime numbers, i.e
$$
 M=pq
$$
where $p$ and $q$ are two prime numbers.
\par To decompose an integer into the product of some prime numbers is
a fundamental problem in number theory([3-4]). How to decompose an
integer into the product of some prime number is difficult problem([5]). There
are several integer decomposition methods, such as short division,
Chinese phase reduction algorithm and so on([6]). However, existing
methods are still difficult to decompose large numbers.
\par In this paper, a new method of an integer decomposition is proposed. The method presented in
this paper is to transform integer decomposition problem into
solving a system of the linear equations.
\section{A New Method of integer decomposition}
Previous integer factorization methods are based on integer division, in this section, we will give an integer decomposition method based on solving algebraic equations on the set $\{0,1\}$.
\par For convenience, the set of the integers $0$ and $1$ is denoted as
$$
\Re=\{0,1\}
$$
\par Let $\alpha\in\Re$, we write
$$
\bar{\alpha}=1-\alpha
$$
which is called the conjugate number of $\alpha$.
\par It is not difficult to verify the following lemma.
\par {\bf Lemma 1.} Let $\alpha, \beta\in\Re$, then for any positive integer $k$
$$
\alpha\beta\in\Re\mbox{}\hspace{16pt}\alpha^{k}=\alpha\in\Re\mbox{}\hspace{16pt}\bar{\alpha}^{k}
=\bar{\alpha}\in\Re\mbox{}\hspace{16pt}|\beta-\alpha|\in\Re
$$
\par Let $M$ be an odd number, its binary representation is as follows
\begin{equation}
M=2^{N}+\sum\limits^{N-1}_{i=1}\gamma_{i}2^{i}+1
\end{equation}
where $\gamma_{i}\in \Re$ $(i=1,2,\cdots,N-1)$. By the binary representation (1) of $M$, we have
$$
2^{N}<1+2^{N}\leq M\leq 2^{N+1}-1<2^{N+1}
$$
which implies that
\begin{equation}
\log_{2}M-1<N<\log_{2}M
\end{equation}
\par Suppose $M$ has
the following binary decomposition
\begin{equation}
M=2^{N}+\sum\limits^{N-1}_{i=1}\gamma_{i}2^{i}+1=\left(\sum\limits^{m}_{i=0}\alpha_{i}2^{i}\right)\left(\sum\limits^{n}_{j=0}\beta_{j}2^{j}\right)
\end{equation}
where $\alpha_{0}=\beta_{0}=\alpha_{m}=\beta_{n}=1$, $m\geq1$, $n\geq1$, $\alpha_{i}\in \Re$
$(i=1,2,\cdots,m-1)$ and $\beta_{j}\in \Re$ $(j=1,2,\cdots,n-1)$.
\par By the binary decomposition (3) of $M$, we obtain
\begin{equation}
2^{m+n}<(2^{m}+1)(2^{n}+1)\leq M\leq(2^{m+1}-1)(2^{n+1}-1)<2^{m+n+2}
\end{equation}
By the inequality (4), we get
$$
2^{m+n}<M\leq 2^{N+1}-1<2^{N+1}\mbox{}\hspace{16pt}2^{m+n+2}>M\geq 2^{N}+1>2^{N}
$$
From the above two inequalities, we find that
$$
m+n<N+1\mbox{}\hspace{16pt}m+n+2>N
$$
That is
$$
N-2<m+n<N+1
$$
Since $m$ and $n$ are integers, we have $m+n=N-1$ or $m+n=N$. Thus, we have the following theorem.
\par {\bf Theorem 1.} Suppose that
$$
M=2^{N}+\sum\limits^{N-1}_{i=1}\gamma_{i}2^{i}+1=\left(\sum\limits^{m}_{i=0}\alpha_{i}2^{i}\right)\left(\sum\limits^{n}_{j=0}\beta_{j}2^{j}\right)
$$
then $m+n=N-1$ or $m+n=N$, where $\alpha_{0}=\beta_{0}=\alpha_{m}=\beta_{n}=1$, $\alpha_{i}\in \Re$
$(i=1,2,\cdots,m-1)$ and $\beta_{j}\in \Re$ $(j=1,2,\cdots,n-1)$.
\par The following two examples show that both
cases in Theorem 1 can occur.
\par {\bf Example 1.} $2^{5}+2^{3}+2^{2}+1=(2+1)(2^{3}+2^{2}+2+1)$.
In this identity, $m=1$, $n=3$ and $N=5$, thus, $m+n=N-1$.
\par {\bf Example 2.} $2^{4}+2^{3}+2+1=(2+1)(2^{3}+1)$.
In this identity, $m=1$, $n=3$ and $N=4$, thus, $m+n=N$.
\par {\bf Corollary 1.} If
$$
2^{N}+2^{k}+\sum\limits^{k-1}_{i=1}\gamma_{i}2^{i}+1=\left(2^{m}+\sum\limits^{m-1}_{i=1}\alpha_{i}2^{i}+1\right)\left(2^{n}+\sum\limits^{n-1}_{j=1}\beta_{j}2^{j}+1\right)
$$
and $k\leq \frac{N}{2}$, then $m+n=N-1$.
\par {\bf Proof.} By the binary decomposition (3) of $M$, we get
\begin{equation}
2^{N}+2^{k}+\sum\limits^{k-1}_{i=1}\gamma_{i}2^{i}+1 \geq
2^{m+n}+2^{m}+2^{n}+1 > 2^{m+n}+2\sqrt{2^{m+n}}
\end{equation}
If $m+n=N$, by (5), we have
\begin{equation}
2^{N}+2^{k}+\sum\limits^{k-1}_{i=1}\gamma_{i}2^{i}+1 >
2^{N}+2^{\frac{N}{2}+1}
\end{equation}
Using $ k\leq \frac{N}{2}$, we find that
\begin{equation}
2^{N}+2^{k}+\sum\limits^{k-1}_{i=1}\gamma_{i}2^{i}+1 \leq
2^{N}+2^{k+1}-1<2^{N}+2^{\frac{N}{2}+1}
\end{equation}
(6) is contradiction with (7). By Theorem 1, we get $m+n=N-1$.
\par Theorem 1 tells us that the number of unknowns in the binary decomposition (3) of $M$ will not exceed $\log_{2}M$. This raises a question: what conditions should $\alpha_{i}\in \Re(i=1,2,\cdots,m-1)$ and $\beta_{j}\in \Re(j=1,2,\cdots,n-1)$ satisfy to make the decomposition (3) to hold?
\par Next, we will derive the
conditions that $\alpha_{i}\in \Re(i=1,2,\cdots,m-1)$ and $\beta_{j}\in \Re(j=1,2,\cdots,n-1)$ should satisfy.
\par By Theorem 1, without loss of generality, we always assume that
$$
1\leq m \leq \left[\frac{N}{2}\right]\mbox{}\hspace{16pt}m\leq n
$$
In the subsequent derivation, we will use the binary representations of
$m$
$$
m=2^{k_{0}}+\sum\limits^{k_{0}-1}_{i=0}\varepsilon_{i}2^{i}
$$
where $k_{0}$ is a positive integer greater than or equal to 1, and $\varepsilon_{i}\in \Re$,
$i=0,1,\cdots,k_{0}-1$.
\par Expanding the right hand side of (3) yields
$$
2^{N}+\sum\limits^{N-1}_{i=1}\gamma_{i}2^{i}+1=\sum\limits^{m+n}_{k=0}\left(\sum\limits_{i+j=k}\alpha_{i}\beta_{j}\right)2^{k}
$$
which can be written as
\begin{equation}
\begin{array}{lll}
2^{N}+\sum\limits^{N-1}_{i=1}\gamma_{i}2^{i}+1
&=&\alpha_{m}\beta_{n}2^{m+n}\\
&+&(\alpha_{m-1}\beta_{n}+\alpha_{m}\beta_{n-1})2^{m+n-1}\\
&+&\sum\limits^{m+n-2}_{j=n+1}(\sum\limits^{m}_{i=j-n}\alpha_{i}\beta_{j-i})2^{j}\\
&+&\sum\limits^{n}_{j=m}(\sum\limits^{m}_{i=0}\alpha_{i}\beta_{j-i})2^{j}\\
&+&\sum\limits^{m-1}_{j=2}(\sum\limits^{j}_{i=0}\alpha_{i}\beta_{j-i})2^{j}\\
&+&(\alpha_{0}\beta_{1}+\alpha_{1}\beta_{0})2\\
&+&\alpha_{0}\beta_{0}
\end{array}
\end{equation}
\par Since $\alpha_{0}=\beta_{0}=1$, it is clear that (8) can be written as
\begin{equation}
\begin{array}{lll}
2^{N-1}+\sum\limits^{N-1}_{i=2}\gamma_{i}2^{i-1}+\gamma_{1}
&=&\alpha_{m}\beta_{n}2^{m+n-1}\\
&+&(\alpha_{m-1}\beta_{n}+\alpha_{m}\beta_{n-1})2^{m+n-2}\\
&+&\sum\limits^{m+n-2}_{j=n+1}(\sum\limits^{m}_{i=j-n}\alpha_{i}\beta_{j-i})2^{j-1}\\
&+&\sum\limits^{n}_{j=m}(\sum\limits^{m}_{i=0}\alpha_{i}\beta_{j-i})2^{j-1}\\
&+&\sum\limits^{m-1}_{j=2}(\sum\limits^{j}_{i=0}\alpha_{i}\beta_{j-i})2^{j-1}\\
&+&(\alpha_{0}\beta_{1}+\alpha_{1}\beta_{0})\\
\end{array}
\end{equation}
By lemma 1, we have
$$
0\leq\alpha_{0}\beta_{1}+\alpha_{1}\beta_{0}\leq 2
$$
Thus, the binary representation of $\alpha_{0}\beta_{1}+\alpha_{1}\beta_{0}$ can be written as
\begin{equation}
\alpha_{0}\beta_{1}+\alpha_{1}\beta_{0}=\tau_{10}+\tau_{11}2
\end{equation}
where $\tau_{10}, \tau_{11}\in \Re$.
\par {\bf Remark 1.} If we add some higher terms of 2 which the degree is greater than or equal to 2 on the right side equation (10), its solution does not change. This means that equation (10) is equivalent to the following equation
$$
\alpha_{0}\beta_{1}+\alpha_{1}\beta_{0}=\tau_{10}+\tau_{11}2+\sum\limits^{l}_{j=2}\tau_{1l}2^{l}
$$
In fact, all the terms in the above equation are greater than or equal to zero, and its left side is less than or equal to 2, so we have to have
$$
\tau_{12}=\tau_{13}=\cdots=\tau_{1l}=0
$$
For convenience, we call $\tau_{1j}(j=2,3,\cdots,l)$ carry.
\par Substituting (10) into (9) yields
\begin{equation}
\begin{array}{lll}
2^{N-1}+\sum\limits^{N-1}_{i=2}\gamma_{i}2^{i-1}+\gamma_{1}
&=&\alpha_{m}\beta_{n}2^{m+n-1}\\
&+&(\alpha_{m-1}\beta_{n}+\alpha_{m}\beta_{n-1})2^{m+n-2}\\
&+&\sum\limits^{m+n-2}_{j=n+1}(\sum\limits^{m}_{i=j-n}\alpha_{i}\beta_{j-i})2^{j-1}\\
&+&\sum\limits^{n}_{j=m}(\sum\limits^{m}_{i=0}\alpha_{i}\beta_{j-i})2^{j-1}\\
&+&\sum\limits^{m-1}_{j=3}(\sum\limits^{j}_{i=0}\alpha_{i}\beta_{j-i})2^{j-1}\\
&+&(\sum\limits^{2}_{i=0}\alpha_{i}\beta_{2-i}+\tau_{11})2\\
&+&\tau_{10}\\
\end{array}
\end{equation}
In the above equation, except for the last term on both sides, all
the other term can be divided by 2, so we have
$\tau_{10}=\gamma_{1}$. Substituting
$\tau_{10}=\gamma_{1}$ into (10) yields
\begin{equation}
\alpha_{0}\beta_{1}+\alpha_{1}\beta_{0}=\gamma_{1}+\tau_{11}2
\end{equation}
\par Using $\tau_{10}=\gamma_{1}$, (11) can be written as
\begin{equation}
\begin{array}{lll}
2^{N-2}+\sum\limits^{N-1}_{i=3}\gamma_{i}2^{i-2}+\gamma_{2}
&=&\alpha_{m}\beta_{n}2^{m+n-2}\\
&+&(\alpha_{m-1}\beta_{n}+\alpha_{m}\beta_{n-1})2^{m+n-3}\\
&+&\sum\limits^{m+n-2}_{j=n+1}(\sum\limits^{m}_{i=j-n}\alpha_{i}\beta_{j-i})2^{j-2}\\
&+&\sum\limits^{n}_{j=m}(\sum\limits^{m}_{i=0}\alpha_{i}\beta_{j-i})2^{j-2}\\
&+&\sum\limits^{m-1}_{j=4}(\sum\limits^{j}_{i=0}\alpha_{i}\beta_{j-i})2^{j-2}\\
&+&(\sum\limits^{3}_{i=0}\alpha_{i}\beta_{3-i})2\\
&+&\sum\limits^{2}_{i=0}\alpha_{i}\beta_{2-i}+\tau_{11}\\
\end{array}
\end{equation}
In the above equation, $\tau_{11}$ is obtained by carrying. Similarly, because of
$0\leq\sum\limits^{2}_{i=0}\alpha_{i}\beta_{2-i}+\tau_{11}\leq 4$, the  binary representation of $\sum\limits^{2}_{i=0}\alpha_{i}\beta_{2-i}+\tau_{11}$ can be expressed as
$$
\sum\limits^{2}_{i=0}\alpha_{i}\beta_{2-i}+\tau_{11}=\gamma_{2}+\tau_{21}2+\tau_{22}2^{2}
$$
where $\tau_{21}, \tau_{22}\in \Re$. For
the sake of unity, we will write the above equation as follows
\begin{equation}
\sum\limits^{2}_{i=0}\alpha_{i}\beta_{2-i}+\sum\limits^{1}_{k=1}\tau_{k(2-k)}=\gamma_{2}+\sum\limits^{2}_{j=1}\tau_{2j}2^{j}
\mbox{}\hspace{16pt}
\end{equation}
Substituting (14) into (13) yields
$$
\begin{array}{lll}
2^{N-2}+\sum\limits^{N-1}_{i=3}\gamma_{i}2^{i-2}+\gamma_{2}
&=&\alpha_{m}\beta_{n}2^{m+n-2}\\
&+&(\alpha_{m-1}\beta_{n}+\alpha_{m}\beta_{n-1})2^{m+n-3}\\
&+&\sum\limits^{m+n-2}_{j=n+1}(\sum\limits^{m}_{i=j-n}\alpha_{i}\beta_{j-i})2^{j-2}\\
&+&\sum\limits^{n}_{j=m}(\sum\limits^{m}_{i=0}\alpha_{i}\beta_{j-i})2^{j-2}\\
&+&\sum\limits^{m-1}_{j=4}(\sum\limits^{j}_{i=0}\alpha_{i}\beta_{j-i})2^{j-2}\\
&+&(\sum\limits^{3}_{i=0}\alpha_{i}\beta_{3-i})2\\
&+&\gamma_{2}+\sum\limits^{2}_{j=1}\tau_{2j}2^{j}\\
\end{array}
$$
which can be rewritten as
\begin{equation}
\begin{array}{lll}
2^{N-3}+\sum\limits^{N-1}_{i=4}\gamma_{i}2^{i-3}+\gamma_{3}
&=&\alpha_{m}\beta_{n}2^{m+n-3}\\
&+&(\alpha_{m-1}\beta_{n}+\alpha_{m}\beta_{n-1})2^{m+n-4}\\
&+&\sum\limits^{m+n-2}_{j=n+1}(\sum\limits^{m}_{i=j-n}\alpha_{i}\beta_{j-i})2^{j-3}\\
&+&\sum\limits^{n}_{j=m}(\sum\limits^{m}_{i=0}\alpha_{i}\beta_{j-i})2^{j-3}\\
&+&\sum\limits^{m-1}_{j=5}(\sum\limits^{j}_{i=0}\alpha_{i}\beta_{j-i})2^{j-3}\\
&+&(\sum\limits^{4}_{i=0}\alpha_{i}\beta_{4-i}+\tau_{22})2\\
&+&\sum\limits^{3}_{i=0}\alpha_{i}\beta_{3-i}+\tau_{21}\\
\end{array}
\end{equation}
Similarly, by Eq. (15), we have
\begin{equation}
\begin{array}{lll}
\sum\limits^{3}_{i=0}\alpha_{i}\beta_{3-i}+\sum\limits^{1}_{k=1}\tau_{(3-k)k}
=\gamma_{3}+\sum\limits^{2}_{j=1}\tau_{3j}2^{j}
\end{array}
\end{equation}
\par {\bf Remark 2.} If equation (10) is replaced by the equation in the Remark 1, then the carry in equation (16) is $\tau_{21}+\tau_{12}$, by Remark 1, we have $\tau_{21}+\tau_{12}=\tau_{21}$. It can be seen from the previous derivation that as long as the left side of equation (10)
remains unchanged, no matter how many higher-order terms are added to the right side of equation (10) and equation (14), the sum of carrying terms in equation (16) will not change. This is true for all the equations derived later.

\par Substituting (16) into (15) yields
$$
\begin{array}{lll}
2^{N-4}+\sum\limits^{N-1}_{i=5}\gamma_{i}2^{i-4}+\gamma_{4}
&=&\alpha_{m}\beta_{n}2^{m+n-4}\\
&+&(\alpha_{m-1}\beta_{n}+\alpha_{m}\beta_{n-1})2^{m+n-5}\\
&+&\sum\limits^{m+n-2}_{j=n+1}(\sum\limits^{m}_{i=j-n}\alpha_{i}\beta_{j-i})2^{j-4}\\
&+&\sum\limits^{n}_{j=m}(\sum\limits^{m}_{i=0}\alpha_{i}\beta_{j-i})2^{j-4}\\
&+&\sum\limits^{m-1}_{j=6}(\sum\limits^{j}_{i=0}\alpha_{i}\beta_{j-i})2^{j-4}\\
&+&(\sum\limits^{5}_{i=0}\alpha_{i}\beta_{5-i}+\tau_{32})2\\
&+&\sum\limits^{4}_{i=0}\alpha_{i}\beta_{4-i}+\tau_{22}+\tau_{31}\\
\end{array}
$$
Continuing the above process, we obtain
\begin{equation}
\left\{
\begin{array}{lr}
\sum\limits^{4}_{i=0}\alpha_{i}\beta_{4-i}+\sum\limits^{2}_{k=1}\tau_{(4-k)k}
=\gamma_{4}+\sum\limits^{2}_{j=1}\tau_{4j}2^{j} & \\
\sum\limits^{5}_{i=0}\alpha_{i}\beta_{5-i}+\sum\limits^{2}_{k=1}\tau_{(5-k)k}
=\gamma_{5}+\sum\limits^{3}_{j=1}\tau_{5j}2^{j} & \\
\sum\limits^{6}_{i=0}\alpha_{i}\beta_{6-i}+\sum\limits^{2}_{k=1}\tau_{(6-k)k}
=\gamma_{6}+\sum\limits^{3}_{j=1}\tau_{6j}2^{j} & \\
\sum\limits^{7}_{i=0}\alpha_{i}\beta_{7-i}+\sum\limits^{2}_{k=1}\gamma_{(7-k)k}
=\gamma_{7}+\sum\limits^{3}_{j=1}\gamma_{7j}2^{j} &
\end{array}
\right.
\end{equation}
\begin{equation}
\left\{
\begin{array}{lr}
\sum\limits^{8}_{i=0}\alpha_{i}\beta_{8-i}+\sum\limits^{3}_{k=1}\tau_{(8-k)k}
=\gamma_{8}+\sum\limits^{3}_{j=1}\tau_{8j}2^{j} & \\
\sum\limits^{9}_{i=0}\alpha_{i}\beta_{9-i}+\sum\limits^{3}_{k=1}\tau_{9-k)k}
=\gamma_{9}+\sum\limits^{3}_{j=1}\tau_{9j}2^{j} & \\
\sum\limits^{10}_{i=0}\alpha_{i}\beta_{10-i}+\sum\limits^{3}_{k=1}\tau_{(10-k)k}
=\gamma_{10}+\sum\limits^{3}_{j=1}\tau_{10j}2^{j} & \\
\sum\limits^{11}_{i=0}\alpha_{i}\beta_{11-i}+\sum\limits^{3}_{k=1}\tau_{(11-k)k}
=\gamma_{11}+\sum\limits^{3}_{j=1}\tau_{11j}2^{j} & \\
\sum\limits^{12}_{i=0}\alpha_{i}\beta_{12-i}+\sum\limits^{3}_{k=1}\tau_{(12-k)k}
=\gamma_{12}+\sum\limits^{4}_{j=1}\tau_{12j}2^{j} & \\
\sum\limits^{13}_{i=0}\alpha_{i}\beta_{13-i}+\sum\limits^{3}_{k=1}\tau_{(13-k)k}
=\gamma_{13}+\sum\limits^{4}_{j=1}\gamma_{13j}2^{j} & \\
\sum\limits^{14}_{i=0}\alpha_{i}\beta_{14-i}+\sum\limits^{3}_{k=1}\tau_{(14-k)k}
=\gamma_{14}+\sum\limits^{4}_{j=1}\gamma_{14j}2^{j} & \\
\sum\limits^{15}_{i=0}\alpha_{i}\beta_{15-i}+\sum\limits^{3}_{k=1}\tau_{(15-k)k}
=\gamma_{15}+\sum\limits^{4}_{j=1}\gamma_{15j}2^{j}
\end{array}
\right.
\end{equation}
In general, for any positive integer $k$ (with $1<k\leq k_{0}-1$), we
have
\begin{equation}
\left\{
\begin{array}{lr}
\sum\limits^{j}_{i=0}\alpha_{i}\beta_{j-i}+d_{jk}=
\gamma_{j}+\sum\limits^{k}_{i=1}\tau_{ji}2^{i}\mbox{}\hspace{16pt} 2^{k}\leq j<2^{k+1}-k-1\\
\sum\limits^{j}_{i=0}\alpha_{i}\beta_{j-i}+d_{jk}=\gamma_{j}+\sum\limits^{k+1}_{i=1}\tau_{ji}2^{i}\mbox{}\hspace{16pt} 2^{k+1}-k-1\leq j<2^{k+1}
\end{array}
\right.
\end{equation}
where $d_{jk}=\sum\limits^{k}_{i=1}\tau_{(j-i)i}$.
\par Noting that $m<2^{k_{0}+1}$, there are three possibilities.
\par {\bf Case I.} $m< 2^{k_{0}+1}-k_{0}-1$. If $2^{k_{0}}\leq j\leq m$, we find that
$$
0\leq\sum\limits^{j}_{i=0}\alpha_{i}\beta_{j-i}+d_{jk_{0}}\leq m+1+k_{0}< 2^{k_{0}+1}-k_{0}-1+1+k_{0}=2^{k_{0}+1}
$$
Thus, the binary representation of $\sum\limits^{j}_{i=0}\alpha_{i}\beta_{j-i}+d_{jk_{0}}$ can be expressed as
\begin{equation}
\sum\limits^{j}_{i=0}\alpha_{i}\beta_{j-i}+d_{jk_{0}}=
\gamma_{j}+\sum\limits^{k_{0}}_{i=1}\tau_{ji}2^{i}
\end{equation}
where $d_{jk_{0}}=\sum\limits^{k_{0}}_{i=1}\tau_{(j-i)i}$, $j=2^{k_{0}}, 2^{k_{0}}+1, ..., m$.
\par The coefficient before $2^{j}$ has $m+1+k_{0}$ terms if $m< j\leq n$, thus, we get
\begin{equation}
\sum\limits^{m}_{i=0}\alpha_{i}\beta_{j-i}+d_{jk_{0}}=
\gamma_{j}+\sum\limits^{k_{0}}_{i=1}\tau_{ji}2^{i}
\end{equation}
\par If $j>n$, the coefficient of $2^{j}$ is less a term than that of $2^{j-1}$, we have
\begin{equation}
\sum\limits^{m}_{i=j-n}\alpha_{i}\beta_{j-i}+d_{jk_{0}}=
\gamma_{j}+\sum\limits^{k_{0}}_{i=1}\tau_{ji}2^{i}
\end{equation}
\par By (12), (14), and (16)-(22), we have the following theorem.
\par {\bf Theorem 2.} Let $m< 2^{k_{0}+1}-k_{0}-1$ and $2\leq k<k_{0}$. If the decomposition
$$
M=2^{N}+\sum\limits^{N-1}_{i=1}\gamma_{i}2^{i}+1=\left(\sum\limits^{m}_{i=0}\alpha_{i}2^{i}\right)\left(\sum\limits^{n}_{j=0}\beta_{j}2^{j}\right)
$$
holds, then $\alpha_{i}(i=1,2,\cdots,m-1)$ and
$\beta_{j}(j=1,2,\cdots,n-1)$ must satisfy the following equations
\begin{equation}
\left\{
\begin{array}{lr}
\beta_{1}+\alpha_{1}=\gamma_{1}+\tau_{11}2\\
\sum\limits^{2}_{i=0}\alpha_{i}\beta_{2-i}+\tau_{11}=\gamma_{2}+\tau_{21}2+\tau_{22}2^{2}\\
\sum\limits^{3}_{i=0}\alpha_{i}\beta_{3-i}+\tau_{21}=\gamma_{3}+\tau_{31}2+\tau_{32}2^{2}\\
\sum\limits^{j}_{i=0}\alpha_{i}\beta_{j-i}+d_{jk}=
\gamma_{j}+\sum\limits^{k}_{i=1}\tau_{ji}2^{i}\mbox{}\hspace{16pt} 2^{k}\leq j<2^{k+1}-k-1\\
\sum\limits^{j}_{i=0}\alpha_{i}\beta_{j-i}+d_{jk}=\gamma_{j}+\sum\limits^{k+1}_{i=1}\tau_{ji}2^{i}\mbox{}\hspace{16pt} 2^{k+1}-k-1\leq j<2^{k+1}\\
\sum\limits^{j}_{i=0}\alpha_{i}\beta_{j-i}+d_{jk_{0}}=
\gamma_{j}+\sum\limits^{k_{0}}_{i=1}\tau_{ji}2^{i}\mbox{}\hspace{16pt}
2^{k_{0}}\leq j\leq m\\
\sum\limits^{m}_{i=0}\alpha_{i}\beta_{j-i}+d_{jk_{0}}=
\gamma_{j}+\sum\limits^{k_{0}}_{i=1}\tau_{ji}2^{i}\mbox{}\hspace{16pt}
m<j\leq n\\
\sum\limits^{m}_{i=j-n}\alpha_{i}\beta_{j-i}+d_{jk_{0}}=
\gamma_{j}+\sum\limits^{k_{0}}_{i=1}\tau_{ji}2^{i}\mbox{}\hspace{16pt}
j>n
\end{array}
\right.
\end{equation}
where $\alpha_{0}=\beta_{0}=\alpha_{m}=\beta_{n}=1$, $\alpha_{i}\in \Re$
$(i=1,2,\cdots,m-1)$ and $\beta_{j}\in \Re$ $(j=1,2,\cdots,n-1)$, $d_{jk}=\sum\limits^{k}_{i=1}\tau_{(j-i)i}$ which is determined by the first $j-1$ equations of Eqs.(23) and is independent of the t-th equation.
\par Eqs.(23) is a little complicated. Can we write it more simply in form? The following lemma shows that it is true.
\par {\bf Lemma 2.} Suppose that the positive integer $M$ has two binary representations
$$
M=\sum\limits^{N}_{i=0}\gamma_{i}2^{i}=\sum\limits^{N+k}_{i=0}\bar{\gamma}_{i}2^{i}
$$
then
$$
\gamma_{i}=\bar{\gamma}_{i}\mbox{}\hspace{8pt}i=0,1,2,\cdots,N
$$
$$
\bar{\gamma}_{N+1}=\bar{\gamma}_{N+2}=\cdots=\bar{\gamma}_{N+k}=0
$$
\par {\bf Proof.} It is clear that
$$
2\left(\sum\limits^{N+k}_{i=1}\bar{\gamma}_{i}2^{i-1}-\sum\limits^{N}_{i=1}\gamma_{i}2^{i-1}\right)=\gamma_{0}-\bar{\gamma}_{0}
$$
Noting that $\left|\gamma_{0}-\bar{\gamma}_{0}\right|\leq 1$, we get $\gamma_{0}=\bar{\gamma}_{0}$, and the above equation can
be rewritten as
$$
2\left(\sum\limits^{N+k}_{i=2}\bar{\gamma}_{i}2^{i-2}-\sum\limits^{N}_{i=2}\gamma_{i}2^{i-2}\right)=\gamma_{1}-\bar{\gamma}_{1}
$$
which implies that
$$
\gamma_{1}=\bar{\gamma}_{1}
$$
Continuing the above process, we have
$$
\gamma_{i}=\bar{\gamma}_{i}\mbox{}\hspace{8pt}i=0,1,2,\cdots,N
$$
and
$$
\sum\limits^{N+k}_{i=N+1}\bar{\gamma}_{i}2^{i-N-1}=0
$$
The final equation shows that
$$
\bar{\gamma}_{N+1}=\bar{\gamma}_{N+2}=\cdots=\bar{\gamma}_{N+k}=0
$$
\par By Lemma 2, Eqs. (23) is equivalent to the following equations
\begin{equation}
\left\{
\begin{array}{lr}
\beta_{1}+\alpha_{1}=\gamma_{1}+\sum\limits^{2m}_{i=1}\tau_{1i}2^{i}\\
\sum\limits^{2}_{i=0}\alpha_{i}\beta_{2-i}+\tau_{11}=\gamma_{2}+\sum\limits^{2m}_{i=1}\tau_{2i}2^{i}\\
\sum\limits^{3}_{i=0}\alpha_{i}\beta_{3-i}+\tau_{21}=\gamma_{3}+\sum\limits^{2m}_{i=1}\tau_{3i}2^{i}\\
\sum\limits^{j}_{i=0}\alpha_{i}\beta_{j-i}+d_{jk}=
\gamma_{j}+\sum\limits^{2m}_{i=1}\tau_{ji}2^{i}\mbox{}\hspace{16pt}2^{k}\leq j<2^{k+1}-k-1\\
\sum\limits^{j}_{i=0}\alpha_{i}\beta_{j-i}+d_{jk}=\gamma_{j}+\sum\limits^{2m}_{i=1}\tau_{ji}2^{i}\mbox{}\hspace{16pt} 2^{k+1}-k-1\leq j<2^{k+1}\\
\sum\limits^{j}_{i=0}\alpha_{i}\beta_{j-i}+d_{jk_{0}}=
\gamma_{j}+\sum\limits^{2m}_{i=1}\tau_{ji}2^{i}\mbox{}\hspace{16pt}
2^{k_{0}}\leq j< m\\
\sum\limits^{m}_{i=0}\alpha_{i}\beta_{j-i}+d_{jk_{0}}=
\gamma_{j}+\sum\limits^{2m}_{i=1}\tau_{ji}2^{i}\mbox{}\hspace{16pt}
m\leq j\leq n\\
\sum\limits^{m}_{i=j-n}\alpha_{i}\beta_{j-i}+d_{jk_{0}}=
\gamma_{j}+\sum\limits^{2m}_{i=1}\tau_{ji}2^{i}\mbox{}\hspace{16pt}
j>n
\end{array}
\right.
\end{equation}
According to the right side of Eqs.(24) and the carry of binary numbers, we can write Eqs.(24) in the following form
$$
\left\{
\begin{array}{lr}
\beta_{1}+\alpha_{1}=\gamma_{1}+\sum\limits^{2m}_{i=1}\tau_{1i}2^{i}\\
\sum\limits^{2}_{i=0}\alpha_{i}\beta_{2-i}+\bar{d}_{21}=\gamma_{2}+\sum\limits^{2m}_{i=1}\tau_{2i}2^{i}\\
\sum\limits^{3}_{i=0}\alpha_{i}\beta_{3-i}+\bar{d}_{32}=\gamma_{3}+\sum\limits^{2m}_{i=1}\tau_{3i}2^{i}\\
\sum\limits^{j}_{i=0}\alpha_{i}\beta_{j-i}+\bar{d}_{j(j-1)}=
\gamma_{j}+\sum\limits^{2m}_{i=1}\tau_{ji}2^{i}\mbox{}\hspace{16pt} 2^{k}\leq j<2^{k+1}-k-1\\
\sum\limits^{j}_{i=0}\alpha_{i}\beta_{j-i}+\bar{d}_{j(j-1)}=\gamma_{j}+\sum\limits^{2m}_{i=1}\tau_{ji}2^{i}\mbox{}\hspace{16pt} 2^{k+1}-k-1\leq j<2^{k+1}\\
\sum\limits^{j}_{i=0}\alpha_{i}\beta_{j-i}+\bar{d}_{j(j-1)}=
\gamma_{j}+\sum\limits^{2m}_{i=1}\tau_{ji}2^{i}\mbox{}\hspace{16pt}
2^{k_{0}}\leq j\leq m\\
\sum\limits^{m}_{i=0}\alpha_{i}\beta_{j-i}+\bar{d}_{jm}=
\gamma_{j}+\sum\limits^{2m}_{i=1}\tau_{ji}2^{i}\mbox{}\hspace{16pt}
m<j\leq n\\
\sum\limits^{m}_{i=j-n}\alpha_{i}\beta_{j-i}+\bar{d}_{jm}=
\gamma_{j}+\sum\limits^{2m}_{i=1}\tau_{ji}2^{i}\mbox{}\hspace{16pt}
j>n
\end{array}
\right.
$$
which can be written in the following concise form
\begin{equation}
\left\{
\begin{array}{lr}
\beta_{1}+\alpha_{1}=\gamma_{1}+\sum\limits^{2m}_{i=1}\tau_{1i}2^{i}\\
\sum\limits^{j}_{i=0}\alpha_{i}\beta_{j-i}+\bar{d}_{j(j-1)}=
\gamma_{j}+\sum\limits^{2m}_{i=1}\tau_{ji}2^{i}\mbox{}\hspace{16pt} 2 \leq j\leq m\\
\sum\limits^{m}_{i=0}\alpha_{i}\beta_{j-i}+\bar{d}_{jm}=
\gamma_{j}+\sum\limits^{2m}_{i=1}\tau_{ji}2^{i}\mbox{}\hspace{16pt}
m<j\leq n\\
\sum\limits^{m}_{i=j-n}\alpha_{i}\beta_{j-i}+\bar{d}_{jm}=
\gamma_{j}+\sum\limits^{2m}_{i=1}\tau_{ji}2^{i}\mbox{}\hspace{16pt}
j>n
\end{array}
\right.
\end{equation}
where
$$
\bar{d}_{j(j-1)}=\sum\limits^{j-1}_{i=1}\tau_{(j-i)i}\mbox{}\hspace{8pt}2\leq j\leq m
$$
and
$$
\bar{d}_{jm}=\sum\limits^{m}_{i=1}\tau_{(j-i)i}\mbox{}\hspace{8pt} j>m
$$
\par It is clear that $\bar{d}_{jk}$ is determined by the first $j-1$ equations of Eqs.(25) and is independent of the j-th equation.
\par In fact, if we start from the first equation $\beta_{1}+\alpha_{1}=\gamma_{1}+\sum\limits^{m}_{i=1}\tau_{1i}2^{i}$ and repeat the process of derving Eqs.(23), we can also get Eqs.(25).
\par We can also prove directly that Eqs.(25) is equivalent to  Eqs.(23). Because $0\leq \beta_{1}+\alpha_{1}\leq 2$, from the fist equation of Eqs.(25), we have
$$
\tau_{12}=\tau_{13}=\cdots=\tau_{1m}=0
$$
Using $\bar{d}_{21}=\tau_{11}$, we obtain
$$
\tau_{23}=\tau_{24}=\cdots=\tau_{2m}=0
$$
Thus, $\bar{d}_{32}=\tau_{21}$, by the third equation of Eqs.(25), we get
$$
\tau_{33}=\tau_{34}=\cdots=\tau_{3m}=0
$$
\par If we continue this process, we can prove that Eqs.(25) is equivalent to  Eqs.(23).
\par Eqs.(25) is particularly simple in form and does not appear $k_{0}$, so it is convenient for later mathematical derivation.
\par It can be seen from Eqs.(23) and  Eqs.(25) that the left side of their first equation plays a decisive role, so long as the left side of the first equation remains unchanged, and the addition of any high-order term to the right side of other equations and their binary carry do not change the solutions of  Eqs.(23) and  Eqs.(25).
\par {\bf Case II.} $m\geq 2^{k_{0}+1}-k_{0}-1$ and $n<2^{k_{0}+1}$.
\par If $2^{k_{0}+1}-k_{0}-1\leq j<m$. $\sum\limits^{j}_{i=0}\alpha_{i}\beta_{j-i}+d_{jk_{0}}$ has $j+1+k_{0}$ terms. Noting that
$$
j+1+k_{0}\geq 2^{k_{0}+1}-k_{0}-1+1+k_{0}=2^{k_{0}+1}
$$
we have
$$
\sum\limits^{j}_{i=0}\alpha_{i}\beta_{j-i}+d_{jk_{0}}=
\gamma_{j}+\sum\limits^{k_{0}+1}_{i=1}\tau_{ji}2^{i}
$$
\par If $m\leq j\leq n$. $\sum\limits^{j}_{i=0}\alpha_{i}\beta_{j-i}+d_{jk_{0}}$ has $m+1+k_{0}$ terms. It is clear that
$$
m+1+k_{0}\leq m+1+\log_{2}m<2m<2\left(2^{k_{0}}-1\right)<2^{k_{0}+2}
$$
Thus, we have
$$
\sum\limits^{j}_{i=0}\alpha_{i}\beta_{j-i}+d_{jk_{0}}=
\gamma_{j}+\sum\limits^{k_{0}+1}_{i=1}\tau_{ji}2^{i}
$$
\par It is clear that the number of terms in $\sum\limits^{m}_{i=j+1-n}\alpha_{i}\beta_{j+1-i}$ is $1$ less than that in $\sum\limits^{m}_{i=j-n}\alpha_{i}\beta_{j-i}$ if $j\geq n$. So we find that
$$
\sum\limits^{m}_{i=j-n}\alpha_{i}\beta_{j-i}+d_{jk_{0}}=
\gamma_{j}+\sum\limits^{k_{0}+1}_{i=1}\tau_{ji}2^{i}\mbox{}\hspace{16pt}
 n< j<2^{k_{0}+1}
$$
$$
\sum\limits^{m}_{i=j-n}\alpha_{i}\beta_{j-i}+d_{j(k_{0}+1)}=
\gamma_{j}+\sum\limits^{k_{0}+1}_{i=1}\tau_{ji}2^{i}\mbox{}\hspace{16pt}
 2^{k_{0}+1}\leq j\leq n+k_{0}+1
$$
\par If $j>n+k_{0}+1$, the number of terms in $\sum\limits^{m}_{i=j-n}\alpha_{i}\beta_{j-i}+d_{j(k_{0}+1)}$ is $m+n-j+k_{0}+2$. It is clear that
$$
m+n-j+k_{0}+2<m+n-(n+k_{0}+1)+k_{0}+2=m+1\leq 2^{k_{0}+1}
$$
which implies that
$$
\sum\limits^{m}_{i=j-n}\alpha_{i}\beta_{j-i}+d_{j(k_{0}+1)}=
\gamma_{j}+\sum\limits^{k_{0}}_{i=1}\tau_{ji}2^{i}\mbox{}\hspace{16pt}
 k_{0}<j-n-1\leq2k_{0}
$$
$$
\sum\limits^{m}_{i=j-n}\alpha_{i}\beta_{j-i}+d_{jk_{0}}=
\gamma_{j}+\sum\limits^{k_{0}}_{i=1}\gamma_{i}^{j}2^{i}\mbox{}\hspace{16pt}
 j-n-1>2k_{0}
$$
\par {\bf Theorem 3.} Let  $m\geq 2^{k_{0}+1}-k_{0}-1$, $n<2^{k_{0}+1}$ and $2\leq k<k_{0}$. If the decomposition
$$
M=2^{N}+\sum\limits^{N-1}_{i=1}\gamma_{i}2^{i}+1=\left(\sum\limits^{m}_{i=0}\alpha_{i}2^{i}\right)\left(\sum\limits^{n}_{j=0}\beta_{j}2^{j}\right)
$$
holds, then $\alpha_{i}(i=1,2,\cdots,m-1)$ and
$\beta_{j}(j=1,2,\cdots,n-1)$ must satisfy the following equations
\begin{equation}
\left\{
\begin{array}{lr}
\beta_{1}+\alpha_{1}=\gamma_{1}+\tau_{11}2\\
\sum\limits^{2}_{i=0}\alpha_{i}\beta_{2-i}+\tau_{11}=\gamma_{2}+\tau_{21}2+\tau_{22}2^{2}\\
\sum\limits^{3}_{i=0}\alpha_{i}\beta_{3-i}+\tau_{21}
=\gamma_{3}+\tau_{31}2+\tau_{32}2^{2}\\
\sum\limits^{j}_{i=0}\alpha_{i}\beta_{j-i}+d_{jk}=
\gamma_{j}+\sum\limits^{k}_{i=1}\tau_{ji}2^{i}\mbox{}\hspace{16pt} 2^{k}\leq j<2^{k+1}-k-1\\
\sum\limits^{j}_{i=0}\alpha_{i}\beta_{j-i}+d_{jk}=\gamma_{j}+\sum\limits^{k+1}_{i=1}\tau_{ji}2^{i}\mbox{}\hspace{16pt} 2^{k+1}-k-1\leq j<2^{k+1}\\
\sum\limits^{j}_{i=0}\alpha_{i}\beta_{j-i}+d_{jk_{0}}=
\gamma_{j}+\sum\limits^{k_{0}}_{i=1}\tau_{ji}2^{i}\mbox{}\hspace{16pt}
2^{k_{0}}\leq j<2^{k_{0}+1}-k_{0}-1\\
\sum\limits^{j}_{i=0}\alpha_{i}\beta_{j-i}+d_{jk_{0}}=
\gamma_{j}+\sum\limits^{k_{0}+1}_{i=1}\tau_{ji}2^{i}\mbox{}\hspace{16pt}
 2^{k_{0}+1}-k_{0}-1\leq j \leq n\\
\sum\limits^{m}_{i=j-n}\alpha_{i}\beta_{j-i}+d_{jk_{0}}=
\gamma_{j}+\sum\limits^{k_{0}+1}_{i=1}\tau_{ji}2^{i}\mbox{}\hspace{16pt}
 n< j<2^{k_{0}+1}\\
\sum\limits^{m}_{i=j-n}\alpha_{i}\beta_{j-i}+d_{j(k_{0}+1)}=
\gamma_{j}+\sum\limits^{k_{0}+1}_{i=1}\tau_{ji}2^{i}\mbox{}\hspace{16pt}
 2^{k_{0}+1}\leq j\leq n+k_{0}+1\\
\sum\limits^{m}_{i=j-n}\alpha_{i}\beta_{j-i}+d_{j(k_{0}+1)}=
\gamma_{j}+\sum\limits^{k_{0}}_{i=1}\tau_{ji}2^{i}\mbox{}\hspace{16pt}
 k_{0}<j-n-1\leq 2k_{0}\\
\sum\limits^{m}_{i=j-n}\alpha_{i}\beta_{j-i}+d_{jk_{0}}=
\gamma_{j}+\sum\limits^{k_{0}}_{i=1}\tau_{ji}2^{i}\mbox{}\hspace{16pt}
 j>n+2k_{0}+1
\end{array}
\right.
\end{equation}
where $\alpha_{0}=\beta_{0}=\alpha_{m}=\beta_{n}=1$, $\alpha_{i}\in \Re$
$(i=1,2,\cdots,m-1)$ and $\beta_{j}\in \Re$ $(j=1,2,\cdots,n-1)$, $d_{jk}=\sum\limits^{k}_{i=1}\tau_{(j-i)i}$.
\par Similar to Case I, it is not difficult to prove that Eqs. (26) is also equivalent to Eqs. (25).

\par {\bf Case III.} $m\geq 2^{k_{0}+1}-k_{0}-1$ and $n\geq2^{k_{0}+1}$.
\par Similar to case II, we have the following theorem.
\par {\bf Theorem 4.} Let  $m\geq 2^{k_{0}+1}-k_{0}-1$, $n\geq2^{k_{0}+1}$ and $2\leq k<k_{0}$. If the decomposition
$$
M=2^{N}+\sum\limits^{N-1}_{i=1}\gamma_{i}2^{i}+1=\left(\sum\limits^{m}_{i=0}\alpha_{i}2^{i}\right)\left(\sum\limits^{n}_{j=0}\beta_{j}2^{j}\right)
$$
holds, then $\alpha_{i}(i=1,2,\cdots,m-1)$ and
$\beta_{j}(j=1,2,\cdots,n-1)$ must satisfy the following equations
\begin{equation}
\left\{
\begin{array}{lr}
\beta_{1}+\alpha_{1}=\gamma_{1}+\tau_{11}2\\
\sum\limits^{2}_{i=0}\alpha_{i}\beta_{2-i}+\tau_{11}=\gamma_{2}+\tau_{21}2+\tau_{22}2^{2}\\
\sum\limits^{3}_{i=0}\alpha_{i}\beta_{3-i}+\tau_{21}
=\gamma_{3}+\tau_{31}2+\tau_{32}2^{2}\\
\sum\limits^{j}_{i=0}\alpha_{i}\beta_{j-i}+d_{jk}=
\gamma_{j}+\sum\limits^{k}_{i=1}\tau_{ji}2^{i}\mbox{}\hspace{16pt} 2^{k}\leq j<2^{k+1}-k-1\\
\sum\limits^{j}_{i=0}\alpha_{i}\beta_{j-i}+d_{jk}=\gamma_{j}+\sum\limits^{k+1}_{i=1}\tau_{ji}2^{i}\mbox{}\hspace{16pt} 2^{k+1}-k-1\leq j<2^{k+1}\\
\sum\limits^{j}_{i=0}\alpha_{i}\beta_{j-i}+d_{jk_{0}}=
\gamma_{j}+\sum\limits^{k_{0}}_{i=1}\tau_{ji}2^{i}\mbox{}\hspace{16pt}
2^{k_{0}}\leq j< 2^{k_{0}+1}-k_{0}-1\\
\sum\limits^{j}_{i=0}\alpha_{i}\beta_{j-i}+d_{jk_{0}}=
\gamma_{j}+\sum\limits^{k_{0}+1}_{i=1}\tau_{ji}2^{i}\mbox{}\hspace{16pt}
 2^{k_{0}+1}-k_{0}-1\leq j < 2^{k_{0}+1}\\
\sum\limits^{j}_{i=0}\alpha_{i}\beta_{j-i}+d_{j(k_{0}+1)}=
\gamma_{j}+\sum\limits^{k_{0}+1}_{i=1}\tau_{ji}2^{i}\mbox{}\hspace{16pt}
 2^{k_{0}+1}\leq j \leq n\\
\sum\limits^{m}_{i=j-n}\alpha_{i}\beta_{j-i}+d_{j(k_{0}+1)}=
\gamma_{j}+\sum\limits^{k_{0}+1}_{i=1}\tau_{ji}2^{i}\mbox{}\hspace{16pt}
n<j\leq n+k_{0}+1\\
\sum\limits^{m}_{i=j-n}\alpha_{i}\beta_{j-i}+d_{j(k_{0}+1)}=
\gamma_{j}+\sum\limits^{k_{0}}_{i=1}\tau_{ji}2^{i}\mbox{}\hspace{16pt}
k_{0}<j-n-1\leq2k_{0}\\
\sum\limits^{m}_{i=j-n}\alpha_{i}\beta_{j-i}+d_{jk_{0}}=
\gamma_{j}+\sum\limits^{k_{0}}_{i=1}\tau_{ji}2^{i}\mbox{}\hspace{16pt}
j>n+2k_{0}+1
\end{array}
\right.
\end{equation}
where $\alpha_{0}=\beta_{0}=\alpha_{m}=\beta_{n}=1$, $\alpha_{i}\in \Re$
$(i=1,2,\cdots,m-1)$ and $\beta_{j}\in \Re$ $(j=1,2,\cdots,n-1)$, $d_{jk}=\sum\limits^{k}_{i=1}\tau_{(j-i)i}$.
\par Similarly, Eqs. (27) is equivalent to Eqs. (25).
\par Eqs. (25) is the unified form of Eqs. (23), Eqs. (26) and Eqs. (27).
\par If Eqs. (25) hold, by Theorem 1, the following corollary holds.
\par {\bf Corollary 2.} If $j>N$, then
$$
d_{j(k_{0}+1)}=\sum\limits^{k_{0}+1}_{i=1}\tau_{(j-i)i}=0
$$
\par It is clear that $d_{j(k_{0}+1)}=0$ if and only if
$$
\tau_{(j-1)1}=\tau_{(j-2)2}=\cdots=\tau_{(j-k_{0}-1)(k_{0}+1)}=0
$$
\par Noting that $m\geq 3$ and $k_{0}\leq \log_{2}m$, we obtain
$$
N-k_{0}-1\geq n+m-\log_{2}m-1
$$
It is not hard to verify that if $m\geq 3$, $f(m)=m-\log_{2}m$ is increasing function, we have
$$
f(m)=m-\log_{2}m\geq f(3)=3-l\log_{2}3>3-\log_{2}4=1
$$
which implies that
$$
N-k_{0}-1\geq n+m-\log_{2}m-1>n
$$
\par The following two corollaries hold.
\par \par {\bf Corollary 3.} If $m+n=N$, then the necessary and sufficient condition for the decomposition formula (3) to hold is that the equations
\begin{equation}
\left\{
\begin{array}{lr}
\beta_{1}+\alpha_{1}=\gamma_{1}+\sum\limits^{2m}_{i=1}\tau_{1i}2^{i}\\
\sum\limits^{j}_{i=0}\alpha_{i}\beta_{j-i}+\bar{d}_{j(j-1)}=
\gamma_{j}+\sum\limits^{2m}_{i=1}\tau_{ji}2^{i}\mbox{}\hspace{16pt} 2 \leq j\leq m\\
\sum\limits^{m}_{i=0}\alpha_{i}\beta_{j-i}+\bar{d}_{jm}=
\gamma_{j}+\sum\limits^{2m}_{i=1}\tau_{ji}2^{i}\mbox{}\hspace{16pt}
m<j\leq n\\
\sum\limits^{m}_{i=j-n}\alpha_{i}\beta_{j-i}+\bar{d}_{jm}=
\gamma_{j}+\sum\limits^{2m}_{i=1}\tau_{ji}2^{i}\mbox{}\hspace{16pt}
j>n\\
\alpha_{N-1}+\beta_{N-1}+d_{(N-1)(k_{0+1})}=\gamma_{N-1}\\
d_{N(k_{0}+1)}=0\\
\end{array}
\right.
\end{equation}
has a solution.
\par \par {\bf Corollary 4.} If $m+n=N-1$, then the necessary and sufficient condition for the decomposition formula (3) to hold is that the equations
\begin{equation}
\left\{
\begin{array}{lr}
\beta_{1}+\alpha_{1}=\gamma_{1}+\sum\limits^{2m}_{i=1}\tau_{1i}2^{i}\\
\sum\limits^{j}_{i=0}\alpha_{i}\beta_{j-i}+\bar{d}_{j(j-1)}=
\gamma_{j}+\sum\limits^{2m}_{i=1}\tau_{ji}2^{i}\mbox{}\hspace{16pt} 2 \leq j\leq m\\
\sum\limits^{m}_{i=0}\alpha_{i}\beta_{j-i}+\bar{d}_{jm}=
\gamma_{j}+\sum\limits^{2m}_{i=1}\tau_{ji}2^{i}\mbox{}\hspace{16pt}
m<j\leq n\\
\sum\limits^{m}_{i=j-n}\alpha_{i}\beta_{j-i}+\bar{d}_{jm}=
\gamma_{j}+\sum\limits^{2m}_{i=1}\tau_{ji}2^{i}\mbox{}\hspace{16pt}
j>n\\
d_{N(k_{0}+1)}=1
\end{array}
\right.
\end{equation}
has a solution.
\par Corollary 3 and Corollary 4 tell us that if we want the decomposition (3) to be true, the unknowns $\alpha_{i}(i=1,2,\cdots,m-1)$ and $\beta_{j}(j=1,2,\cdots,n-1)$ must satisfy the system of multivariate quadratic equations (28) or (29) on the set $\{0,1\}$. Thus, the integer decomposition problem is transformed into solving a system of the algebraic equations.
\par By Corollary 3 and Corollary 4, we have the following two remarks.
\par {\bf Remark 3.} The necessary and sufficient condition for an odd number $M$ to be a composite number is that there exists $m$
and $n$ with $m+n=N$ such that Eqs.(28) has a solution or $m$ and $n$ with $m+n=N-1$ such that Eqs.(29) has a solution.
\par {\bf Remark 4.} The necessary and sufficient condition for an odd number $M$ to be a prime number if and only if Eqs.(28) has no solution for any $m$ and $n$ with $m+n=N$ and Eqs.(29) has no solution for $m$ and $n$ with $m+n=N-1$.
\section{Some Important Lemmas}
\par In the last section, we established the system of multivariate quadratic algebraic Eqs.(28) and Eqs.(29) for the unknowns $\alpha_{i}(i=1,2,\cdots,m-1)$ and $\beta_{j}(j=1,2,\cdots,n-1)$ on the set $\Re$. It is clear that the existing methods for solving quadratic algebraic equations cannot used to solve them, so we have to find new method.
\par In this section, we will give several important lemmas, which are the basis of the method that we will propose for solving Eqs.(28) and Eqs.(29). These lemmas ensure that the operations are performed on the set $\Re$.
\par{\bf Lemma 3.} Let $\alpha, \beta\in\Re$, then
\par{\bf (1)} $\alpha+\beta=|\beta-\alpha|+\alpha\beta2$
\par{\bf (2)} $\left|\alpha-\beta\right|\beta=\bar{\alpha}\beta$
\par{\bf (3)} $\alpha\bar{\alpha}=0$\mbox{}\hspace{8pt}$\bar{\bar{\alpha}}=\alpha$\mbox{}\hspace{8pt}$|\alpha-\bar{\alpha}|=|1-2\alpha|=1$
\par{\bf (4)} $\alpha\beta\left|\alpha-\beta\right|=0$
\par{\bf (5)} $\left|\left|\beta-\alpha\right|-\alpha\right|=\beta$
\par{\bf (6)} If $f: \Re\longmapsto \Re$, then
$$
f(\alpha)=(1-\alpha)f(0)+\alpha f(1)=\bar{\alpha}f(0)+\alpha f(1)
$$
\par{\bf (7)} If $f: \Re\times \Re\longmapsto \Re$, then
$$
f(\alpha,\beta)=\bar{\alpha}\bar{\beta}f(0,0)+\bar{\alpha}\beta f(0,1)+\alpha\bar{\beta}f(1,0)+\alpha\beta f(1,1)
$$
\par{\bf (8)} $|\alpha\beta-|\beta-\alpha||=\bar{\alpha}\beta+\alpha=|\bar{\alpha}\beta-\alpha|$
\par{\bf Proof.} {\bf Proof of (1).} By Lemma 1, it is clear that
$$
\alpha+\beta=\alpha^{2}+\beta^{2}=(\beta-\alpha)^{2}+\alpha\beta2=|\beta-\alpha|+\alpha\beta2
$$
\par {\bf Proof of (2).} It is clear that
$$
\left|\alpha-\beta\right|\beta=\left|\alpha\beta-\beta^{2}\right|=\left|\alpha\beta-\beta\right|=(1-\alpha)\beta=\bar{\alpha}\beta
$$
\par {\bf Proof of (5).} It is not hard to verift that
$$
\left|\left|\beta-\alpha\right|-\alpha\right|=
\left\{
\begin{array}{lr}
\left|\beta\right|=\beta\mbox{}\hspace{16pt} \alpha=0\\
\left|\left|\beta-1\right|-1\right|=\left|1-\beta-1\right|=\beta\mbox{}\hspace{16pt} \alpha=1\\
\end{array}
\right.
$$
\par {\bf Proof of (7).} We computer
$$
\begin{array}{lll}
f(\alpha,\beta)
&=&\bar{\alpha}f(0,\beta)+\alpha f(1,\beta)\\
&=&\bar{\alpha}(\bar{\beta}f(0,0)+\beta f(0,1))+\alpha(\bar{\beta}f(1,0)+\beta f(1,1))\\
&=&\bar{\alpha}\bar{\beta}f(0,0)+\bar{\alpha}\beta f(0,1)+\alpha\bar{\beta}f(1,0)+\alpha\beta f(1,1)\\
\end{array}
$$
\par {\bf Proof of (8).} Writing $f(\alpha,\beta)=|\alpha\beta-|\beta-\alpha||$, by (7), we get
$$
f(\alpha,\beta)=\bar{\alpha}\beta+\alpha\bar{\beta}+\alpha\beta =\bar{\alpha}\beta+\alpha(\bar{\beta}+\beta)=\bar{\alpha}\beta+\alpha
$$
Noting that $\bar{\alpha}\alpha=0$, we have
$$
\bar{\alpha}\beta+\alpha=(\bar{\alpha}\beta)^{2}+\alpha^{2}-2\bar{\alpha}\beta\alpha=(\bar{\alpha}\beta-\alpha)^{2}=\left|\bar{\alpha}\beta-\alpha\right|
$$
\par {\bf Lemma 4.} Let $\beta, \gamma_{1}, \gamma_{2}\in\Re$, if $\left|\beta-\gamma_{1}\right|=\left(\beta-\gamma_{1}\right)^{2}=\gamma_{2}$, then
$$
\beta=(\gamma_{2}-\gamma_{1})^{2}=|\gamma_{2}-\gamma_{1}|
$$
\par {\bf Proof.} If $\gamma_{1}=0$, $(\beta-\gamma_{1})^{2}=\gamma_{2}$ becomes $\beta^{2}=\gamma_{2}$, thus $\beta=\beta^{2}=\gamma_{2}$. If $\gamma_{1}=1$, $(\beta-\gamma_{1})^{2}=\gamma_{2}$ becomes $(1-\beta)^{2}=\gamma_{2}$, thus $1-\beta=\gamma_{2}$,  $\beta=1-\gamma_{2}$.
Combine the two situations, by Lemma 3,  we get
$$
\beta=(1-\gamma_{1})\gamma_{2}+\gamma_{1}(1-\gamma_{2})=\gamma_{1}+\gamma_{2}-2\gamma_{1}\gamma_{2}=(\gamma_{2}-\gamma_{1})^{2}=|\gamma_{2}-\gamma_{1}|
$$
\par For the sake of simplicity, we will introduce some symbols. Let $\delta_{i} \in\Re(i=1,2,\cdots,n)$, for any positive integers $i$ and $j$ with $1\leq i<j\leq n$, define
$$
l_{ii}(\delta_{h})=\delta_{i}
$$
$$
l_{ij}(\delta_{h})=l_{ij}(\delta_{i},\delta_{i+1},\cdots,\delta_{j})=|\delta_{j}-|\delta_{j-1}-\cdots-|\delta_{i+2}-|\delta_{i+1}-\delta_{i}
\underbrace{||\cdots|}_{(j-i)\mbox{\scriptsize times}}
$$
By Lemma 1, it is clear that $l_{ij}(\delta_{h})\in\Re$. We find that
$$
l_{i(i+1)}(\delta_{h})=l_{i(i+1)}(\delta_{i},\delta_{i+1})=|\delta_{i+1}-\delta_{i}|=|\delta_{i+1}-l_{ii}(\delta_{h})|
$$
$$
l_{i(i+2)}(\delta_{h})=l_{i(i+2)}(\delta_{i},\delta_{i+1},\delta_{i+2})=|\delta_{i+2}-|\delta_{i+1}-\delta_{i}||=|\delta_{i+2}-l_{i(i+1)}(\delta_{h})|
$$
$$
l_{i(i+3)}(\delta_{h})=|\delta_{i+3}-|\delta_{i+2}-|\delta_{i+1}-\delta_{i}|||=|\delta_{i+3}-l_{i(i+2)}(\delta_{h})|
$$
In a general way, for any positive integer $k$ with $i+1\leq k\leq n$, we have
$$
l_{ik}(\delta_{h})=|\delta_{k}-l_{i(k-1)}(\delta_{h})|
$$
\par Let $\delta_{1i}\in\Re(i=1,2,\cdots,n)$, write
$$
\delta_{2i}=\delta_{1i}l_{1(i-1)}(\delta_{1h})\mbox{}\hspace{16pt}i=2,3,\cdots,n
$$
$$
\delta_{3i}=\delta_{2i}l_{2(i-1)}(\delta_{2h})\mbox{}\hspace{16pt}i=3,4,\cdots,n
$$
$$
\delta_{4i}=\delta_{3i}l_{3(i-1)}(\delta_{3h})\mbox{}\hspace{16pt}i=4,5,\cdots,n
$$
Popularly, for any positive integer $k$ with $k\leq n$, define
$$
\delta_{ki}=\delta_{(k-1)i}l_{(k-1)(i-1)}(\delta_{(k-1)h})\mbox{}\hspace{16pt}k=2,3,\cdots,n\mbox{}\hspace{16pt}i=k,k+1,\cdots,n
$$
By Lemma 1, we have $\delta_{ki}\in\Re$.
\par We will see that the key to solve the Eqs.(28) and Eqs.(29) is how to express the sum of several numbers on $\Re$ into binary form.
The following lemma provides a method to express  the sum of $n$ numbers on  $\Re$ into binary form.
\par {\bf Lemma 5.} Let $\delta_{1i}\in\Re(i=1,2,\cdots,n)$ , then
\begin{equation}
\sum\limits^{n}_{i=1}\delta_{1i}=\sum\limits^{n}_{i=1}l_{in}(\delta_{ih})2^{i-1}
\end{equation}
\par {\bf Proof.} By Lemma 3, we have
\begin{equation}
\begin{array}{lll}
\sum\limits^{n}_{i=1}\delta_{1i}&=&\delta_{11}+\delta_{12}+\sum\limits^{n}_{i=3}\delta_{1i}\\
&=&\left|\delta_{12}-\delta_{11}\right|+\delta_{12}\delta_{11}2+\sum\limits^{n}_{i=3}\delta_{1i}\\
&=&l_{12}(\delta_{1h})+\delta_{12}l_{11}(\delta_{1h})2+\sum\limits^{n}_{i=3}\delta_{1i}\\
&=&l_{12}(\delta_{1h})+\delta_{13}+\delta_{12}l_{11}(\delta_{1h})2+\sum\limits^{n}_{i=4}\delta_{1i}\\
&=&\left|\delta_{13}-l_{12}(\delta_{1h})\right|+\left(\delta_{12}l_{11}(\delta_{1h})+\delta_{13}l_{12}(\delta_{1h})\right)2+\sum\limits^{n}_{i=4}\delta_{1i}\\
&=&l_{13}(\delta_{1h})+\left(\delta_{12}l_{11}(\delta_{1h})+\delta_{13}l_{12}(\delta_{1h})\right)2+\sum\limits^{n}_{i=4}\delta_{1i}\\
&=&l_{13}(\delta_{1h})+\delta_{14}+\left(\delta_{12}l_{11}(\delta_{1h})+\delta_{13}l_{12}(\delta_{1h})\right)2+\sum\limits^{n}_{i=5}\delta_{1i}\\
&=&l_{14}(\delta_{1h})+\left(\sum\limits^{4}_{i=2}\delta_{1i}l_{1(i-1)}(\delta_{1h})\right)2+\sum\limits^{n}_{i=5}\delta_{1i}\\
&=&\cdots\cdots\cdots\cdots\cdots\cdots\cdots\cdots\cdots\cdots\cdots\cdots\cdots\cdots\cdots\\
&=&l_{1(n-1)}(\delta_{1h})+\delta_{1n}+\left(\sum\limits^{n-1}_{i=2}\delta_{1i}l_{1(i-1)}(\delta_{1h})\right)2\\
&=&l_{1n}(\delta_{1h})+\left(\sum\limits^{n}_{i=2}\delta_{1i}l_{1(i-1)}(\delta_{1h})\right)2\\
\end{array}
\end{equation}
By Eq. (31), Eq. (30) becomes
\begin{equation}
l_{1n}(\delta_{1h})+\left(\sum\limits^{n}_{i=2}\delta_{1i}l_{1(i-1)}(\delta_{1h})\right)2=c_{1}+\sum\limits^{n}_{i=2}c_{i}2^{i-1}
\end{equation}
which implies $2\mid (l_{1n}(\delta_{1h})-c_{1})$. Noting that $|l_{1n}(\delta_{1h})-c_{1}|\leq 1$, we have
$$
c_{1}=l_{1n}(\delta_{1h})
$$
and Eq. (32) can be written as
$$
\sum\limits^{n}_{i=2}\delta_{1i}l_{1(i-1)}(\delta_{1h})=c_{2}+\sum\limits^{n}_{i=3}c_{i}2^{i-2}
$$
That is
$$
\sum\limits^{n}_{i=2}\delta_{2i}=c_{2}+\sum\limits^{n}_{i=3}c_{i}2^{i-2}
$$
Similarly, we get
$$
l_{2n}(\delta_{2h})+\left(\sum\limits^{n}_{i=3}\delta_{2i}l_{2(i-1)}(\delta_{2h})\right)2=c_{2}+\sum\limits^{n}_{i=3}c_{i}2^{i-2}
$$
Thus
$$
c_{2}=l_{2n}(\delta_{2h})
$$
and
\begin{equation}
\sum\limits^{n}_{i=3}\delta_{2i}l_{2(i-1)}(\delta_{2h})=c_{3}+\sum\limits^{n}_{i=4}c_{i}2^{i-3}
\end{equation}
Eq. (33) can be written as
\begin{equation}
\sum\limits^{n}_{i=3}\delta_{3i}=c_{3}+\sum\limits^{n}_{i=4}c_{i}2^{i-3}
\end{equation}
which implies that
$$
c_{3}=l_{3n}(\delta_{3h})
$$
and
$$
\sum\limits^{n}_{i=4}\delta_{4i}=c_{4}+\sum\limits^{n}_{i=5}c_{i}2^{i-4}
$$
Similarly, for any positive integer $k$ with $k\leq n-2$, we have
$$
c_{k}=l_{kn}(\delta_{kh})
$$
and
$$
\sum\limits^{n}_{i=k+1}\delta_{(k+1)i}=c_{k+1}+\sum\limits^{n}_{i=k+2}c_{i}2^{i-k-1}
$$
In particular, if $k=n-2$, we get
$$
c_{n-2}=l_{(n-2)n}(\delta_{(n-2)h})
$$
and
$$
\delta_{(n-1)(n-1)}+\delta_{(n-1)n}=c_{n-1}+c_{n}2
$$
which implies that
$$
c_{n-1}=l_{(n-1)n}(\delta_{(n-1)h})
$$
$$
c_{n}=\delta_{(n-1)n}\delta_{(n-1)(n-1)}=\delta_{(n-1)n}l_{(n-1)(n-1))}(\delta_{(n-1)h})=\delta_{nn}=l_{nn}(\delta_{nh})
$$
\par It is clear that $1+1=0+2$. In this case, $c_{n}=c_{2}=1$. If $n>2$, we computer
$$
\begin{array}{lll}
c_{n}&=&\delta_{(n-1)n}\delta_{(n-1)(n-1)}\\
&=&\delta_{(n-2)n}l_{(n-2)(n-1)}(\delta_{(n-2)h})\delta_{(n-2)(n-1)}l_{(n-2)(n-2)}(\delta_{(n-2)h})\\
&=&\delta_{(n-2)n}\delta_{(n-2)(n-1)}l_{(n-2)(n-1)}(\delta_{(n-2)h})l_{(n-2)(n-2)}(\delta_{(n-2)h})
\end{array}
$$
By Lemma 3, we have
$$
\begin{array}{lll}
& &\delta_{(n-2)(n-1)}l_{(n-2)(n-1)}(\delta_{(n-2)h})l_{(n-2)(n-2)}(\delta_{(n-2)h})\\
&=&\delta_{(n-2)(n-1)}|\delta_{(n-2)(n-1)}-l_{(n-2)(n-2)}(\delta_{(n-2)h})|
$$l_{(n-2)(n-2)}(\delta_{(n-2)h})\\
&=&0
\end{array}
$$
which implies that if $ n>2$, $c_{n}=0$.
\par In fact, if $2^{k}\leq\sum\limits^{n}_{i=1}\delta_{1i}<2^{k+1}$, then for any $j$ with $k<j\leq n$, $c_{j}=0$.
\par {\bf Lemma 6.} Let $\delta_{1i}\in\Re(i=1,2,\cdots,n)$ and $\gamma\in\Re$ be known, then for $\beta, \mu_{i}\in\Re$, the solution of the following equation
\begin{equation}
\beta+\sum\limits^{n}_{i=1}\delta_{1i}=
\gamma+\sum\limits^{n}_{i=1}\mu_{i}2^{i}
\end{equation}
is given by
$$
\beta=\left|\gamma-l_{1n}(\delta_{1h})\right|
$$
$$
\mu_{k}=\left|l_{(k+1)n}(\delta_{(k+1)h})-\bar{\gamma}\prod\limits^{k}_{i=1} l_{in}(\delta_{ih})\right|\mbox{}\hspace{8pt}k=1,2,\cdots,n-1
$$
$$
\mu_{n}=\bar{\gamma}\prod\limits^{n}_{i=1} l_{in}(\delta_{ih})
$$
\par {\bf Proof.} By Lemma 5, Eq. (35) can be written as
\begin{equation}
\beta+l_{1n}(\delta_{1h})+\sum\limits^{n}_{i=2}l_{in}(\delta_{ih})2^{i-1}=
\gamma+\sum\limits^{n}_{i=1}\mu_{i}2^{i}
\end{equation}
By Lemma 3, Eq. (36) can be written as
\begin{equation}
(\beta-l_{1n}(\delta_{1h}))^{2}+\beta l_{1n}(\delta_{1h})2+\sum\limits^{n}_{i=2}l_{in}(\delta_{ih})2^{i-1}=
\gamma+\sum\limits^{n}_{i=1}\mu_{i}2^{i}
\end{equation}
which implies that
$$
(\beta-l_{1n}(\delta_{1h}))^{2}=\gamma
$$
By Lemma 4, we have
\begin{equation}
\beta=\left|\gamma-l_{1n}(\delta_{1h})\right|
\end{equation}
Substituting (38) into (37) yields
\begin{equation}
\left|\gamma-l_{1n}(\delta_{1h})\right| l_{1n}(\delta_{1h})+l_{2n}(\delta_{2h})+\sum\limits^{n}_{i=3}l_{in}(\delta_{ih})2^{i-2}=\mu_{1}+
\sum\limits^{n}_{i=2}\mu_{i}2^{i-1}
\end{equation}
By Lemma 3, Eq. (39) can be expresses as
$$
(1-\gamma)l_{1n}(\delta_{1h})+l_{2n}(\delta_{2h})+\sum\limits^{n}_{i=3}l_{in}(\delta_{ih})2^{i-2}=\mu_{1}+
\sum\limits^{n}_{i=2}\mu_{i}2^{i-1}
$$
That is
$$
\bar{\gamma}l_{1n}(\delta_{1h})+l_{2n}(\delta_{2h})+\sum\limits^{n}_{i=3}l_{in}(\delta_{ih})2^{i-2}=\mu_{1}+
\sum\limits^{n}_{i=2}\mu_{i}2^{i-1}
$$
which implies that
$$
\left|l_{2n}(\delta_{2h})-\bar{\gamma}l_{1n}(\delta_{1h})\right|+\bar{\gamma}\prod\limits^{2}_{i=1} l_{in}(\delta_{ih})2
+\sum\limits^{n}_{i=3}l_{in}(\delta_{ih})2^{i-2}=\mu_{1}+
\sum\limits^{n}_{i=2}\mu_{i}2^{i-1}
$$
Thus, we get
$$
\mu_{1}=\left|l_{2n}(\delta_{2h})-\bar{\gamma}l_{1n}(\delta_{1h})\right|
$$
and
\begin{equation}
\bar{\gamma}\prod\limits^{2}_{i=1} l_{in}(\delta_{ih})+l_{3n}(\delta_{3h})
+\sum\limits^{n}_{i=4}l_{in}(\delta_{ih})2^{i-3}=\mu_{2}+
\sum\limits^{n}_{i=3}\gamma_{i}2^{i-2}
\end{equation}
By Eq. (40), we have
$$
\left|l_{3n}(\delta_{3h})-\bar{\gamma}\prod\limits^{2}_{i=1} l_{in}(\delta_{ih})\right|+\bar{\gamma}\prod\limits^{3}_{i=1} l_{in}(\delta_{ih})2
+\sum\limits^{n}_{i=4}l_{in}(\delta_{ih})2^{i-3}=\mu_{2}+
\sum\limits^{n}_{i=3}\mu_{i}2^{i-2}
$$
which implies that
$$
\mu_{2}=\left|l_{3n}(\delta_{3h})-\bar{\gamma}\prod\limits^{2}_{i=1} l_{in}(\delta_{ih})\right|
$$
and
$$
\bar{\gamma}\prod\limits^{3}_{i=1} l_{in}(\delta_{ih})+l_{4n}(\delta_{4h})
+\sum\limits^{n}_{i=5}l_{in}(\delta_{ih})2^{i-4}=\mu_{3}+
\sum\limits^{n}_{i=4}\mu_{i}2^{i-3}
$$
Similarly, we have
$$
\mu_{3}=\left|l_{4n}(\delta_{4h})-\bar{\gamma}\prod\limits^{3}_{i=1} l_{in}(\delta_{ih})\right|
$$
and
$$
\bar{\gamma}\prod\limits^{4}_{i=1} l_{in}(\delta_{ih})+l_{5n}(\delta_{5h})
+\sum\limits^{n}_{i=6}l_{in}(\delta_{ih})2^{i-5}=\mu_{4}+
\sum\limits^{n}_{i=5}\mu_{i}2^{i-4}
$$
Continuing the process,  for any positive integer $k$ with $k\leq n-3$, we obtain
$$
\mu_{k}=\left|l_{(k+1)n}(\delta_{(k+1)h})-\bar{\gamma}\prod\limits^{k}_{i=1} l_{in}(\delta_{ih})\right|
$$
and
$$
\bar{\gamma}\prod\limits^{k+1}_{i=1} l_{in}(\delta_{ih})+l_{(k+2)n}(\delta_{(k+2)h})
+\sum\limits^{n}_{i=k+3}l_{in}(\delta_{ih})2^{i-k-2}=\mu_{k+1}+
\sum\limits^{n}_{i=k+2}\mu_{i}2^{i-k-1}
$$
If $k=n-3$, we have
$$
\bar{\gamma}\prod\limits^{n-2}_{i=1} l_{in}(\delta_{ih})+l_{(n-1)n}(\delta_{(n-1)h})
+l_{nn}(\delta_{nh})2=\mu_{n-2}+\mu_{n-1}2+\mu_{n}2^{2}
$$
By the above equation, we obtain
$$
\mu_{n-2}=\left|l_{(n-1)n}(\delta_{(n-1)h})-\bar{\gamma}\prod\limits^{n-2}_{i=1} l_{in}(\delta_{ih})\right|
$$
and
$$
\bar{\gamma}\prod\limits^{n-1}_{i=1} l_{in}(\delta_{ih})+l_{nn}(\delta_{nh})=\mu_{n-1}+\mu_{n}2
$$
Thus, we have
$$
\mu_{n-1}=\left|l_{nn}(\delta_{nh})-\bar{\gamma}\prod\limits^{n-1}_{i=1} l_{in}(\delta_{ih})\right|
$$
and
$$
\mu_{n}=\bar{\gamma}\prod\limits^{n}_{i=1} l_{in}(\delta_{ih})
$$
\par Because $l_{12}(\delta_{1h})=\left|\delta_{12}-\delta_{11}\right|=(\delta_{12}-\delta_{11})^{2}$, computing  $l_{12}(\delta_{1h})$ requires a subtraction and a multiplication, which means that it takes two operations. Similarly,  computing  $l_{13}(\delta_{1h})=\left|\delta_{13}-l_{12}(\delta_{1h})\right|$ requires two operations, computing  $l_{14}(\delta_{1h})=\left|\delta_{14}-l_{13}(\delta_{1h})\right|$ requires two operations, continuing the process, computing  $l_{1n}(\delta_{1h})=\left|\delta_{1n}-l_{1(n-1)}(\delta_{1h})\right|$ requires two operations. Thus, the total number of operations
to calculate $l_{1n}(\delta_{1h})$ is $2(n-1)$.
\par The total number of operations for $\delta_{22}=\delta_{12}l_{11}(\delta_{1h})$, $\delta_{23}=\delta_{13}l_{12}(\delta_{1h})$, $\cdots$, and $\delta_{2n}=\delta_{1n}l_{1(n-1)}(\delta_{1h})$ is $(n-1)$ . Thus, computing  $l_{2n}(\delta_{2h})$ requires $n-1+2(n-2)$ operations.
\par Similarly, computing  $l_{3n}(\delta_{3h})$ requires $n-2+2(n-3)$ operations.
\par In general, for any positive integer $i$ (with $2\leq i\leq n$), computing  $l_{in}(\delta_{ih})$ requires $n-i+1+2(n-i)$ operations.
\par To sum up, the total number of operations for $l_{1n}(\delta_{1h})$, $l_{2n}(\delta_{2h})$, $l_{3n}(\delta_{3h})$, $\cdots$, and $l_{nn}(\delta_{nh})$ is
$$
\begin{array}{lll}
& &2(n-1)+[(n-1)+2(n-2)]+[(n-2)+2(n-3)]+\cdots\\
&+&\{n-(n-1)+1+2[n-(n-1)]\}+[n-n+1+2(n-n)]\\
&=&3(n-1)+3(n-2)+3(n-3)+\cdots+3[n-(n-2)]+3[n-(n-1)]\\
&=&\frac{3(n-1)n}{2}
\end{array}
$$
\par Let
$$\rho_{k}=\bar{\gamma}\prod\limits^{k}_{i=1} l_{in}(\delta_{ih})(k=1,2,\cdots,n)$$
It is clear that
$$
\rho_{k+1}=l_{(k+1)n}(\delta_{(k+1)h})\rho_{k}
$$
\par It takes an operation to calculate $\rho_{k}(1\leq k\leq n)$, two for $\beta$ and $\mu_{k}(1\leq n-1)$, and one for $\mu_{n}$. The total number of calculations is
$$
n+2+2(n-1)+1=3n+1
$$
\par {\bf Corollary 5.} The number of calculations required to solve Eq. (35) is
$$
1+\frac{3n(n+1)}{2}
$$
\par {\bf Lemma 7.} Let $\delta_{1i}\in\Re(i=1,2,\cdots,n)$ and $\gamma\in\Re$ be known. Then, a necessary and sufficient condition for the equation
\begin{equation}
\sum\limits^{n}_{i=1}\delta_{1i}=
\gamma+\sum\limits^{n-1}_{i=1}\mu_{i}2^{i}
\end{equation}
to have a solution is
$$
g(\gamma)=\gamma-l_{1n}(\delta_{1h})=0
$$
and if the condition is true, the solution is given by
$$
\mu_{i}=l_{(i+1)n}(\delta_{(i+1)h})\mbox{}\hspace{16pt}i=1,2,\cdots,n-1
$$
\par {\bf Proof.} By Lemma 5, Eq. (41) can be written as
$$
l_{1n}(\delta_{1h})+\sum\limits^{n}_{i=2}l_{in}(\delta_{ih})2^{i-1}=\gamma+\sum\limits^{n}_{i=1}\mu_{i}2^{i}
$$
The necessary and sufficient condition for the above equation to have a solution is $2|(\gamma-l_{1n}(\delta_{1h}))$. Noting that
$|\gamma-l_{1n}(\delta_{1h})|\leq 1$, thus, the necessary and sufficient condition for $2|(\gamma-l_{1n}(\delta_{1h}))$ is $\gamma=l_{1n}(\delta_{1h})$.
\par {\bf Example 3.}  Find out the condition that an odd number
$$
M=2^{N}+\sum\limits^{N-1}_{i=1}\gamma_{i}2^{i}+1
$$
cab be divisible by 17, where $\gamma_{i}\in \Re$
$(i=1,2,\cdots,N-1)$ and $N>8$.
\par Noting that $17=2^{4}+1$ , if $M$ can be divisible by 17, we have
$$
M=2^{N}+\sum\limits^{N-1}_{i=1}\gamma_{i}2^{i}+1=\left(2^{4}+1\right)\left(1+\sum\limits^{n-1}_{j=1}\beta_{j}2^{j}+2^{n}\right)
$$
Noting that $m=4=2^{2}<5=2^{3}-2-1$, by the Theorem 2, $\beta_{j}(j=1,2,\cdots,n-1)$ must satisfy the following equations
\begin{equation}
\left\{
\begin{array}{lr}
\beta_{1}=\gamma_{1}\\
\beta_{2}=\gamma_{2}\\
\beta_{3}=\gamma_{3}\\
\beta_{4}+1=\gamma_{4}+\tau_{41}2\\
\beta_{5}+\beta_{1}+\tau_{41}=\gamma_{5}+\tau_{51}2\\
\beta_{6}+\beta_{2}+\tau_{51}=\gamma_{6}+\tau_{61}2\\
\cdots\cdots\cdots\cdots\cdots\cdots\cdots\cdots\cdots\\
\beta_{n-1}+\beta_{n-5}+\tau_{(n-2)1}=\gamma_{n-1}+\tau_{(n-1)1}2\\
1+\beta_{n-4}+\tau_{(n-1)1}=\gamma_{n}+\tau_{n1}2\\
\beta_{n-3}+\tau_{n1}=\gamma_{n+1}+\tau_{(n+1)1}2\\
\beta_{n-2}+\tau_{(n+1)1}=\gamma_{n+2}+\tau_{(n+2)1}2\\
\beta_{n-1}+\tau_{(n+2)1}=\gamma_{n+3}+\tau_{(n+3)1}2\\
1+\tau_{(n+3)1}=\gamma_{n+4}+\tau_{(n+4)1}2\\
\end{array}
\right.
\end{equation}
\par Eqs.(42) is an over-determined system of equations, we can uniquely solve $\beta_{1}, \beta_{2}, \cdots, \beta_{n-1}$ from its first $n-1$ equations. By Lemma 6, we get
\begin{equation}
\left\{
\begin{array}{lr}
\beta_{1}=\gamma_{1}\\
\beta_{2}=\gamma_{2}\\
\beta_{3}=\gamma_{3}\\
\beta_{4}=\tau_{41}=1-\gamma_{4}=\bar{\gamma}_{4}\\
\beta_{5}=|\gamma_{5}-|\tau_{41}-\beta_{1}||=|\gamma_{5}-|\bar{\gamma}_{4}-\gamma_{1}||\\
\tau_{51}=|\tau_{41}\beta_{1}-\bar{\gamma}_{5}|\tau_{41}-\beta_{1}||=|\bar{\gamma}_{4}\gamma_{1}-\bar{\gamma}_{5}|\bar{\gamma}_{4}-\gamma_{1}||\\
\beta_{6}=|\gamma_{6}-|\tau_{51}-\beta_{2}||\\
\tau_{61}=|\tau_{51}\beta_{2}-\bar{\gamma}_{6}|\tau_{51}-\beta_{2}||\\
\beta_{7}=|\gamma_{7}-|\tau_{61}-\beta_{3}||\\
\tau_{71}=|\tau_{61}\beta_{3}-\bar{\gamma}_{7}|\tau_{61}-\beta_{3}||\\
\cdots\cdots\cdots\cdots\cdots\cdots\cdots\cdots\cdots\\
\beta_{n-1}=|\gamma_{n-1}-|\tau_{(n-2)1}-\beta_{n-5}||\\
\tau_{(n-1)1}=|\tau_{(n-2)1}\beta_{n-5}-\bar{\gamma}_{n-1}|\tau_{(n-2)1}-\beta_{n-5}||\\
\end{array}
\right.
\end{equation}
\par If $n+4=N-1$, by Lemma 7 and Corollary 4, the necessary and sufficient condition for $17|M$ is as follows
\begin{equation}
\left\{
\begin{array}{lr}
g(\gamma_{n})=\gamma_{n}-|\tau_{(n-1)1}-\bar{\beta}_{n-4}|=0\\
g(\gamma_{n+1})=\gamma_{n+1}-|\tau_{n1}-\beta_{n-3}|=0\\
g(\gamma_{n+2})=\gamma_{n+2}-|\tau_{(n+1)1}-\beta_{n-2}|=0\\
g(\gamma_{n+3})=\gamma_{n+3}-|\tau_{(n+2)1}-\beta_{n-1}|=0\\
g(\gamma_{n+4})=1-|1-\tau_{(n+3)1}|\\
\tau_{(n+4)1}=1\\
\end{array}
\right.
\end{equation}
and if (44) hold, we obtain
$$
\left\{
\begin{array}{lr}
\tau_{n1}=|\beta_{n-4}-\bar{\beta}_{n-4}\tau_{(n-1)1}|\\
\tau_{(n+1)1}=\tau_{n1}\beta_{n-3}\\
\tau_{(n+2)1}=\tau_{(n+1)1}\beta_{n-2}\\
\tau_{(n+3)1}=\tau_{(n+2)1}\beta_{n-1}\\
\tau_{(n+4)1}=\tau_{(n+3)1}\\
\end{array}
\right.
$$
\par If $n+4=N$, then $\gamma_{n+4}=1$. By Lemma 7 and Corollary 3, the necessary and sufficient condition for $17|M$ is as follows
\begin{equation}
\left\{
\begin{array}{lr}
g(\gamma_{n})=\gamma_{n}-|\tau_{(n-1)1}-\bar{\beta}_{n-4}|=0\\
g(\gamma_{n+1})=\gamma_{n+1}-|\tau_{n1}-\beta_{n-3}|=0\\
g(\gamma_{n+2})=\gamma_{n+2}-|\tau_{(n+1)1}-\beta_{n-2}|=0\\
g(\gamma_{n+3})=\gamma_{n+3}-|\tau_{(n+2)1}-\beta_{n-1}|=0\\
g(\gamma_{n+4})=1-|1-\tau_{(n+3)1}|=0\\
\tau_{(n+4)1}=0\\
\end{array}
\right.
\end{equation}
and if (45) hold, we obtain
$$
\left\{
\begin{array}{lr}
\tau_{n1}=|\beta_{n-4}-\bar{\beta}_{n-4}\tau_{(n-1)1}|\\
\tau_{(n+1)1}=\tau_{n1}\beta_{n-3}\\
\tau_{(n+2)1}=\tau_{(n+1)1}\beta_{n-2}\\
\tau_{(n+3)1}=\tau_{(n+2)1}\beta_{n-1}\\
\tau_{(n+4)1}=\tau_{(n+3)1}\\
\end{array}
\right.
$$
\par We consider two special cases.
\par {\bf Case A.} $\gamma_{1}=\gamma_{2}=\cdots=\gamma_{N-1}=0$. By Theorem 1, we have $N+4=N-1$
\par {\bf (1).} If $n-1=4(2p_{0}+1)+q_{0}$ with $0\leq q_{0}\leq 3$, by Eqs.(43), we get
$$
\beta_{1}=\beta_{2}=\beta_{3}=0
$$
and for $0\leq p \leq p_{0}-1$
$$
\beta_{4(2p+1)}=\beta_{4(2p+1)+1}=\beta_{4(2p+1)+2}=\beta_{4(2p+1)+3}=1
$$
$$
\tau_{[4(2p+1)]1}=\tau_{[4(2p+1)+1]1}=\tau_{[4(2p+1)+2]1}=\tau_{[4(2p+1)+3]1}=1
$$
$$
\beta_{4(2p+2)}=\beta_{4(2p+2)+1}=\beta_{4(2p+2)+2}=\beta_{4(2p+2)+3}=0
$$
$$
\tau_{[4(2p+2)]1}=\tau_{[4(2p+2)+1]1}=\tau_{[4(2p+2)+3]1}=\tau_{[4(2p+2)+3]1}=1
$$
and
$$
\beta_{4(2p_{0}+1)}=\cdots=\beta_{4(2p_{0}+1)+q_{0}}=1
$$
$$
\tau_{[4(2p_{0}+1)]1}=\cdots=\tau_{[4(2p_{0}+1)+q_{0}]1}=1
$$
\par {\bf (A1).} If $q_{0}=0$, then $\beta_{n-4}=\beta_{n-3}=\beta_{n-2}=0$. By (44), we get
$$
g(\gamma_{n})=\gamma_{n}-|\tau_{(n-1)1}-\bar{\beta}_{n-4}|=0-|1-1|=0
$$
which implies that $\tau_{n1}=1$. We computer
$$
g(\gamma_{n+1})=\gamma_{n+1}-|\tau_{n1}-\beta_{n-3}|=0-|1-0|=1
$$
In this case, Eqs.(42) have no solution.
\par {\bf (A2).} If $q_{0}=1$, then $\beta_{n-4}=\beta_{n-3}=0, \beta_{n-2}=1$.By (44), we obtain
$$
g(\gamma_{n})=0\mbox{}\hspace{16pt}\tau_{n1}=1\mbox{}\hspace{16pt}\gamma_{n+1}=1
$$
In this case, Eqs.(42) have no solution.
\par {\bf (A3).} If $q_{0}=2$, then $\beta_{n-5}=\beta_{n-4}=0,\beta_{n-3}=\beta_{n-2}=1$. By (44), we have
$$
\left\{
\begin{array}{lr}
\gamma_{n}=0\mbox{}\hspace{16pt}\tau_{n1}=1\\
\gamma_{n+1}=0\mbox{}\hspace{16pt}\tau_{(n+1)1}=1\\
\gamma_{n+2}=0\mbox{}\hspace{16pt}\tau_{(n+2)1}=1\\
\gamma_{n+3}=0\mbox{}\hspace{16pt}\tau_{(n+3)1}=1\\
\gamma_{n+4}=0\mbox{}\hspace{16pt}\tau_{(n+4)1}=1\\
\end{array}
\right.
$$
which implies that $d_{N1}=\tau_{(n+4)1}=1$. Thus, in this case, Eqs.(42) have a solution.
\par It is clear that
$$
n=4(2p_{0}+1)+3
$$
$$
N=n+5=4(2p_{0}+1)+8=4[2(p_{0}+1)+1]
$$
\par {\bf (A4).} If $q_{0}=3$, then $\beta_{n-4}=\beta_{n-3}=\beta_{n-2}=1$. By (44), we have
$$
g(\gamma_{n})=\gamma_{n}-|\tau_{(n-1)1}-\bar{\beta}_{n-4}|=0-|1-0|=1
$$
In this case, Eqs.(42) have no solution.
\par {\bf (2).} If $n-1=4\cdot 2p_{0}+q_{0}$ with $0\leq q_{0}\leq 3$, by Eqs.(43), we get
$$
\beta_{1}=\beta_{2}=\beta_{3}=0
$$
and for $0\leq p \leq p_{0}-2$
$$
\beta_{4(2p+1)}=\beta_{4(2p+1)+1}=\beta_{4(2p+1)+2}=\beta_{4(2p+1)+3}=1
$$
$$
\tau_{[4(2p+1)]1}=\tau_{[4(2p+1)+1]1}=\tau_{[4(2p+1)+2]1}=\tau_{[4(2p+1)+3]1}=1
$$
$$
\beta_{4(2p+2)}=\beta_{4(2p+2)+1}=\beta_{4(2p+2)+2}=\beta_{4(2p+2)+3}=0
$$
$$
\tau_{[4(2p+2)]1}=\tau_{[4(2p+2)+1]1}=\tau_{[4(2p+2)+3]1}=\tau_{[4(2p+2)+3]1}=1
$$
and
$$
\beta_{4(2p_{0}-1)}=\beta_{4(2p_{0}-1)+1}=\beta_{4(2p_{0}-1)+2}=\beta_{4(2p_{0}-1)+3}=1
$$
$$
\tau_{4(2p_{0}-1)1}=\tau_{[4(2p_{0}-1)+1]1}=\tau_{[4(2p_{0}-1)+2]1}=\tau_{[4(2p_{0}-1)+3]1}=1
$$
$$
\beta_{4\cdot2p_{0}}=\beta_{4\cdot2p_{0}+1}=\cdots=\beta_{4\cdot2p_{0}+q_{0}}=0
$$
$$
\tau_{(4\cdot2p_{0})1}=\tau_{(4\cdot2p_{0}+1)1}=\cdots=\tau_{(4\cdot2p_{0}+q_{0})1}=1
$$
\par {\bf (A5).} If $q_{0}=0$, then $\beta_{n-4}=\beta_{n-3}=\beta_{n-2}=1$. By (44), we obtain
$$
g(\gamma_{n})=\gamma_{n}-|\tau_{(n-1)1}-\bar{\beta}_{n-4}|=0-|1-0|=1
$$
In this case, Eqs.(42) have no solution.
\par {\bf (A6).} If $q_{0}=1$, then $\beta_{n-4}=\beta_{n-3}=1$, $\beta_{n-2}=0$. By (44), we obtain
$$
g(\gamma_{n})=1
$$
In this case, Eqs.(42) have no solution.
\par {\bf (A7).} If $q_{0}=2$, then $\beta_{n-4}=1$, $\beta_{n-3}=\beta_{n-2}=0$. By (44), we obtain
$$
g(\gamma_{n})=1
$$
In this case, Eqs.(42) have no solution.
\par {\bf (A8).} If $q_{0}=3$, then $\beta_{n-4}=\beta_{n-3}=\beta_{n-2}=0$. By (44), we obtain
$$
g(\gamma_{n})=0\mbox{}\hspace{16pt}\tau_{n1}=1\mbox{}\hspace{16pt}
$$
which implies that
$$
g(\gamma_{n+1})=\gamma_{n+1}-|\tau_{n1}-\beta_{n-3}|=0-|1-0|=1
$$
In this case, Eqs.(42) have no solution.
\par By (A1)-(A8), $2^{N}+1$ can be divisible by 17 if and only if $N=4(2k+1)$, where $k$ is an arbitrary positive integer.
\par {\bf Case B.} $\gamma_{1}=\gamma_{2}=\cdots=\gamma_{N-1}=1$.
\par {\bf (3).} If $n-1=4(2p_{0}+1)+q_{0}$ with $0\leq q_{0}\leq 3$, by Eqs.(43), we get
$$
\beta_{1}=\beta_{2}=\beta_{3}=1
$$
and for $0\leq p \leq p_{0}-1$
$$
\beta_{4(2p+1)}=\beta_{[4(2p+1)+1]}=\beta_{[4(2p+1)+2]}=\beta_{[4(2p+1)+3]}=0
$$
$$
\tau_{4(2p+1)1}=\tau_{[4(2p+1)+1]1}=\tau_{[4(2p+1)+2]1}=\tau_{[4(2p+1)+3]1}=0
$$
$$
\beta_{4(2p+2)}=\beta_{[4(2p+2)+1]}=\beta_{[4(2p+2)+2]}=\beta_{[4(2p+2)+3]}=1
$$
$$
\tau_{4(2p+2)1}=\tau_{[4(2p+2)+1]1}=\tau_{[4(2p+2)+2]1}=\tau_{[4(2p+2)+3]1}=0
$$
and
$$
\beta_{4(2p_{0}+1)}=\beta_{[4(2p_{0}+1)+1]}=\cdots=\beta_{[4(2p_{0}+1)+q_{0}]}=0
$$
$$
\tau_{4(2p_{0}+1)1}=\tau_{[4(2p_{0}+1)+1]1}=\cdots=\tau_{[4(2p_{0}+1)+q_{0}]1}=0
$$
\par {\bf (B1).} If $q_{0}=0$, then $\beta_{n-4}=\beta_{n-3}=\beta_{n-2}=1$. We obtain
$$
g(\gamma_{n})=1
$$
In this case, Eqs.(42) have no solution.
\par {\bf (B2).} If $q_{0}=1$, then $\beta_{n-4}=1=\beta_{n-3}=1$, $\beta_{n-2}=0$. We get
$$
g(\gamma_{n})=1
$$
In this case, Eqs.(42) have no solution.
\par {\bf (B3).} If $q_{0}=2$, then $\beta_{n-4}=1$, $\beta_{n-3}=\beta_{n-2}=0$. We find that
$$
g(\gamma_{n})=1
$$
In this case, Eqs.(42) have no solution.
\par {\bf (B4).} If $q_{0}=3$, then $\beta_{n-4}=\beta_{n-3}=\beta_{n-2}=0$. We have
$$
g(\gamma_{n})=0\mbox{}\hspace{16pt}\tau_{n1}=0\mbox{}\hspace{16pt}g(\gamma_{n+1})=1
$$
In this case, Eqs.(42) have no solution.
\par {\bf (4).} If $n-1=4\cdot 2p_{0}+q_{0}$ with $0\leq q_{0}\leq 3$, by Eqs.(43), we get
$$
\beta_{1}=\beta_{2}=\beta_{3}=1
$$
and for $0\leq p \leq p_{0}-2$
$$
\beta_{4(2p+1)}=\beta_{4(2p+1)+1}=\beta_{4(2p+1)+2}=\beta_{4(2p+1)+3}=0
$$
$$
\tau_{[4(2p+1)]1}=\tau_{[4(2p+1)+1]1}=\tau_{[4(2p+1)+2]1}=\tau_{[4(2p+1)+3]1}=0
$$
$$
\beta_{4(2p+2)}=\beta_{4(2p+2)+1}=\beta_{4(2p+2)+2}=\beta_{4(2p+2)+3}=1
$$
$$
\tau_{[4(2p+2)]1}=\tau_{[4(2p+2)+1]1}=\tau_{[4(2p+2)+3]1}=\tau_{[4(2p+2)+3]1}=0
$$
and
$$
\beta_{4(2p_{0}-1)}=\beta_{4(2p_{0}-1)+1}=\beta_{4(2p_{0}-1)+2}=\beta_{4(2p_{0}-1)+3}=0
$$
$$
\tau_{4(2p_{0}-1)1}=\tau_{[4(2p_{0}-1)+1]1}=\tau_{[4(2p_{0}-1)+2]1}=\tau_{[4(2p_{0}-1)+3]1}=0
$$
$$
\beta_{4\cdot2p_{0}}=\beta_{4\cdot2p_{0}+1}=\cdots=\beta_{4\cdot2p_{0}+q_{0}}=1
$$
$$
\tau_{(4\cdot2p_{0})1}=\tau_{(4\cdot2p_{0}+1)1}=\cdots=\tau_{(4\cdot2p_{0}+q_{0})1}=0
$$
\par {\bf (B5).} If $q_{0}=0$, then $\beta_{n-4}=\beta_{n-3}=\beta_{n-2}=0$. We have
$$
g(\gamma_{n})=0\mbox{}\hspace{16pt}\tau_{n1}=0\mbox{}\hspace{16pt}g(\gamma_{n+1})=1
$$
In this case, Eqs.(42) have no solution.
\par {\bf (B6).} If $q_{0}=1$, then $\beta_{n-4}=\beta_{n-3}=0$, $\beta_{n-2}=1$. We have
$$
g(\gamma_{n})=0\mbox{}\hspace{16pt}\tau_{n1}=0\mbox{}\hspace{16pt}g(\gamma_{n+1})=1
$$
In this case, Eqs.(42) have no solution.
\par {\bf (B7).} If $q_{0}=2$, then $\beta_{n-4}=0$, $\beta_{n-3}=\beta_{n-2}=1$. We have
$$
\begin{array}{lll}
g(\gamma_{n})=0\mbox{}\hspace{16pt}\tau_{n1}=0\\
g(\gamma_{n+1})=0\mbox{}\hspace{16pt}\tau_{(n+1)1}=0\\
g(\gamma_{n+2})=0\mbox{}\hspace{16pt}\tau_{(n+2)1}=0\\
g(\gamma_{n+3})=0\mbox{}\hspace{16pt}\tau_{(n+3)1}=0\\
g(\gamma_{n+4})=0\mbox{}\hspace{16pt}\tau_{(n+4)1}=0\\
\end{array}
$$
\par There are two situations.
\par {\bf Case (B71).} If $n+4=N-1$, by (44),  Eqs.(42) have no solution.
\par {\bf Case (B72).} If $n+4=N$, by (45),  Eqs.(42) have solution. It is not clear that
$$
N=n+4=4\cdot 2p_{0}+7=8(p_{0}+1)-1
$$
\par {\bf (B8).} If $q_{0}=3$, then $\beta_{n-4}=\beta_{n-3}=\beta_{n-2}=1$.
We have
$$
g(\gamma_{n})=1
$$
In this case, Eqs.(42) have no solution.
\par Noting that
$$
2^{N}+2^{N=1}+\cdots+1=2^{N+1}-1
$$
By (B1)-(B8),  $2^{N+1}-1$ can be divisible by 17 if and only if $N=8k-1$, where $k$ is an arbitrary positive integer.
\par It is not hard to verify directly by he traditional method
$$
\begin{array}{lll}
2^{4(2k+1)}+1&=&(16)^{(2k+1)}+1=(17-1)^{(2k+1)}+1\\
&=&17\bar{M}+(-1)^{(2k+1)}+1=17\bar{M}
\end{array}
$$
and
$$
\begin{array}{lll}
2^{8k-1+1}-1&=&(16)^{2k}-1=(17-1)^{2k}-1\\
&=&17\hat{M}+(-1)^{2k}-1=17\hat{M}
\end{array}
$$
which show that the method in this paper is correct. However, compared with the traditional method, the method seems to be much more complicated, which is not the case. It is difficult to find other integer numbers divided by $17$ except $2^{4(2k+1)}+1$ and $2^{8k}-1$ by the traditional method. Theoretically speaking, using (44) and (45) in example 3, all integer numbers that can be divisible by $17$ can be found out by the method in this paper.
\par It can be seen from Eqs.(25) that if one factor in (3) is given, we can judge the existence of its another factor by Lemma 6 and Lemma 7, and if there is, we can solve it.
\section{How To Solve Equations}
\par Since Eqs.(25) is a system of multivariate quadratic algebraic equations on the set $\{0,1\}$, the traditional method of solving algebraic equations cannot be used to solve it.
\par From example 3, we can see that if a factor of (3) is determined, Eqs.(25) becomes a linear system of equations. Next, based on Lemma 6 and Lemma 7, we will propose a method to solve Eqs.(25).
Since $\alpha_{i}\in\Re(i=0,1,2,\cdots,m)$ and $\alpha_{0}=\alpha_{m}=1$, without loss of generality, we might as well assume that
$$
\alpha_{0}=\alpha_{1}=\alpha_{2}=\alpha_{3}=\cdots=\alpha_{k_{1}}=1
$$
$$
\alpha_{k_{1}+1}=\alpha_{k_{1}+2}=\alpha_{k_{1}+3}=\cdots=\alpha_{k_{2}}=0
$$
$$
 \alpha_{k_{2}+1}=\alpha_{k_{2}+2}=\alpha_{k_{2}+3}=\cdots=\alpha_{k_{3}}=1
$$
$$
 \alpha_{k_{3}+1}=\alpha_{k_{3}+2}=\alpha_{k_{3}+3}=\cdots=\alpha_{k_{4}}=0
$$
$$
 \alpha_{k_{4}+1}=\alpha_{k_{4}+2}=\alpha_{k_{4}+3}=\cdots=\alpha_{k_{5}}=1
$$
$$
\alpha_{k_{5}+1}=\alpha_{k_{5}+2}=\alpha_{k_{5}+3}=\cdots=\alpha_{k_{6}}=0
$$
$$
\cdots\cdots\cdots\cdots\cdots\cdots\cdots\cdots\cdots\cdots\cdots\cdots
$$
$$
\alpha_{k_{2q-1}+1}=\alpha_{k_{2q-1}+2}=\cdots=\alpha_{k_{2q}}=0
$$
$$
\alpha_{k_{2q}+1}=\alpha_{k_{2q}+2}=\cdots=\alpha_{m}=1
$$
where $0\leq \alpha_{k_{1}},\alpha_{k_{2}},\cdots,\alpha_{k_{2q}}\leq m-1$.
\par We can arrange $\alpha_{0}$, $\alpha_{1}$, $\cdots$, and $\alpha_{m}$ into the following sequence
\begin{equation}
\begin{array}{lll}
& &\alpha_{0},\alpha_{1},\alpha_{2},\cdots,\alpha_{m}\\
&=&\overbrace{1,1,\cdots,1}^{k_{1}+1},\overbrace{0,0,\cdots,0}^{k_{2}-k_{1}}\overbrace{1,1,\cdots,1}^{k_{3}-k_{2}},\cdots,
\overbrace{0,0,\cdots,0}^{k_{2q}-k_{2q-1}},\overbrace{1,1,\cdots,1}^{m-k_{2q}}
\end{array}
\end{equation}
The number of all such sequence is $2^{m-1}$(Since $\alpha_{0}=\alpha_{m}=1$). Especially if $k_{1}=m-1$, we get
$$
\alpha_{0},\alpha_{1},\alpha_{2},\cdots,\alpha_{m}=\overbrace{1,1,\cdots,1}^{m+1}
$$
and if $k_{1}=0$ and $k_{2}=m-1$, we have
$$
\alpha_{0},\alpha_{1},\alpha_{2},\cdots,\alpha_{m}=1,\overbrace{0,0\cdots,0}^{m-1},1
$$
\par Substituting (46) into Eqs.(25) yields
\begin{equation}
\left\{
\begin{array}{lr}
\beta_{1}+1=\gamma_{1}+\sum\limits^{2m}_{i=1}\tau_{1i}2^{i}\\
\beta_{j}+\sum\limits^{j}_{i=1}\beta_{i-1}+\bar{d}_{j(j-1)}=
\gamma_{j}+\sum\limits^{2m}_{i=1}\tau_{ji}2^{i}\mbox{}\hspace{16pt} 1 < j\leq k_{1}\\
\beta_{j}+\Sigma(1,j)+\bar{d}_{j(j-1)}=
\gamma_{j}+\sum\limits^{2m}_{i=1}\tau_{ji}2^{i}\mbox{}\hspace{16pt} k_{1} < j\leq k_{2}\\
\beta_{j}+\Sigma(1,j)+\sum\limits^{j-k_{2}}_{i=1}\beta_{i-1}+\bar{d}_{j(j-1)}=
\gamma_{j}+\sum\limits^{2m}_{i=1}\tau_{ji}2^{i}\mbox{}\hspace{16pt} k_{2} < j\leq k_{3}\\
\beta_{j}+\Sigma(2,j)+\bar{d}_{j(j-1)}=
\gamma_{j}+\sum\limits^{2m}_{i=1}\tau_{ji}2^{i}\mbox{}\hspace{16pt} k_{3} < j\leq k_{4}\\
\beta_{j}+\Sigma(2,j)+\sum\limits^{j-k_{4}}_{i=1}\beta_{i-1}+\bar{d}_{j(j-1)}=
\gamma_{j}+\sum\limits^{2m}_{i=1}\tau_{ji}2^{i}\mbox{}\hspace{16pt} k_{4} < j\leq k_{5}\\
\beta_{j}+\Sigma(3,j)+\bar{d}_{j(j-1)}=
\gamma_{j}+\sum\limits^{2m}_{i=1}\tau_{ji}2^{i}\mbox{}\hspace{16pt} k_{5} < j\leq k_{6}\\
\cdots\cdots\cdots\cdots\cdots\cdots\cdots\cdots\cdots\cdots\cdots\cdots\cdots\cdots\cdots\cdots\cdots\cdots\cdots\\
\beta_{j}+\Sigma(q,j)+\bar{d}_{j(j-1)}=
\gamma_{j}+\sum\limits^{2m}_{i=1}\tau_{ji}2^{i}\mbox{}\hspace{16pt} k_{2q-1} < j\leq k_{2q}\\
\beta_{j}+\Sigma(q,j)+
\sum\limits^{j-k_{2q}}_{i=1}\beta_{i-1}+\bar{d}_{j(j-1)}=
\gamma_{j}+\sum\limits^{2m}_{i=1}\tau_{ji}2^{i}\mbox{}\hspace{8pt} k_{2q} < j\leq m\\
\beta_{j}+\Sigma(q,j)+\sum\limits^{j-k_{2q}}_{i=j-m+1}\beta_{i-1}+\bar{d}_{j(j-1)}=
\gamma_{j}+\sum\limits^{2m}_{i=1}\tau_{ji}2^{i}\mbox{}\hspace{8pt} m < j\leq n\\
\sum\limits^{k_{1}}_{i=j-n}\beta_{j-i}+\hat{\Sigma}(1,j)+\bar{d}_{jm}=
\gamma_{j}+\sum\limits^{2m}_{i=1}\tau_{ji}2^{i}\mbox{}\hspace{16pt}
n<j\leq n+k_{1}\\
\hat{\Sigma}(1,j)+\bar{d}_{jm}=
\gamma_{j}+\sum\limits^{2m}_{i=1}\tau_{ji}2^{i}\mbox{}\hspace{16pt}
n+k_{1}<j\leq n+k_{2}\\
\sum\limits^{k_{3}}_{i=j-n}\beta_{j-i}+\hat{\Sigma}(2,j)+\bar{d}_{jm}=
\gamma_{j}+\sum\limits^{2m}_{i=1}\tau_{ji}2^{i}\mbox{}\hspace{16pt}
k_{2}<j-n\leq k_{3}\\
\hat{\Sigma}(2,j)+\bar{d}_{jm}=
\gamma_{j}+\sum\limits^{2m}_{i=1}\tau_{ji}2^{i}\mbox{}\hspace{16pt}
k_{3}<j-n\leq k_{4}\\
\cdots\cdots\cdots\cdots\cdots\cdots\cdots\cdots\cdots\cdots\cdots\cdots\cdots\cdots\cdots\\
\hat{\Sigma}(q,j)+\bar{d}_{jm}=
\gamma_{j}+\sum\limits^{2m}_{i=1}\tau_{ji}2^{i}\mbox{}\hspace{16pt}
k_{2q-1}<j-n\leq k_{2q}\\
\sum\limits^{m}_{i=j-n}\beta_{j-i}+\bar{d}_{jm}=
\gamma_{j}+\sum\limits^{2m}_{i=1}\tau_{ji}2^{i}\mbox{}\hspace{16pt}
k_{2q}<j-n\leq m
\end{array}
\right.
\end{equation}
where $\beta_{0}=\beta_{n}=1$, and
$$
\Sigma(k,j)=\sum\limits^{k-1}_{h=0}\sum\limits^{j-k_{2h}}_{i=j-k_{2h+1}+1}\beta_{i-1}\mbox{}\hspace{16pt} k_{0}=0
$$
and
$$
\hat{\Sigma}(l,j)=\sum\limits^{q-1}_{h=l}\sum\limits^{k_{2h+1}}_{i=k_{2h}+1}\beta_{j-i}+\sum\limits^{m}_{i=k_{2q}+1}\beta_{j-i}
$$
\par Therefore, we transform to solve the multivariate quadratic algebraic equations (25) into to solve linear algebraic equations (47). By Lemma 6, we can solve $\beta_{1}$, $\beta_{2}$, $\cdots$ and $\beta_{n-1}$ from the first $n-1$ equations of Eqs.(47). Next, we can use Lemma 7 to judge the existence of solutions of Eqs.(47). If there is no solution, our calculation ends. If there is a solution, we can judge whether the decomposition formula (3) holds by corollary 3 and corollary 4.
\par {\bf Example 4.} Decompose the integer $F_{4}=2^{2^{4}}+1=2^{16}+1$.
\par By Corollary 1, we have $m+n=15$ with $m\leq n$. Seven sets
of solutions are obtain, they are $(m,n)=(1,14)$, $(m,n)=(2,13)$,
$(m,n)=(3,12)$, $(m,n)=(4,11)$, $(m,n)=(5,10)$, $(m,n)=(6,9)$ and
$(m,n)=(7,8)$. Seven cases need to be considered in the example.
\par {\bf Case 1.} $(m,n)=(1,14)$, we have
$$2^{16}+1=(2+1)(2^{14}+\sum\limits^{13}_{i=1}\beta_{i}2^{i}+1)$$
Expanding the right hand side of the above equation yields
$$2^{16}+1=2^{15}+(1+\beta_{13})2^{14}+\sum\limits^{13}_{i=2}(\beta_{i}+\beta_{i-1})2^{i}+(1+\beta_{1})2+1$$
We obtain
\begin{equation}
\left\{
\begin{array}{lll}
1+\beta_{1}=\tau_{11}2\\
\beta_{2}+\beta_{1}+\tau_{11}=\tau_{21}2+\tau_{22}2^{2}\\
\beta_{3}+\beta_{2}+\tau_{21}=\tau_{31}2+\tau_{32}2^{2}\\
\beta_{i}+\beta_{i-1}+\tau_{(i-1)1}+\tau_{(i-2)2}=\tau_{i1}2+\tau_{i2}2^{2} \mbox{}\hspace{12pt} 3\leq i\leq 13\\
1+\beta_{13}+\tau_{(13)1}+\tau_{(12)2}=\tau_{(14)1}2+\tau_{(14)2}2^{2}\\
1+\tau_{(14)1}+\tau_{(13)2}=\tau_{(15)1}2\\
\tau_{(15)1}+\tau_{(14)2}=1
\end{array}
\right.
\end{equation}
By Lemma 6, we have
$$
\left\{
\begin{array}{lll}
\beta_{1}=\tau_{11}=1\mbox{}\hspace{8pt}\beta_{2}=\tau_{22}=0\mbox{}\hspace{8pt}\tau_{21}=1\\
\beta_{3}=\tau_{31}=1\mbox{}\hspace{8pt}\tau_{32}=0\mbox{}\hspace{8pt}\tau_{41}=1\mbox{}\hspace{8pt}\beta_{4}=\tau_{42}=0\\
\beta_{5}=\tau_{51}=1\mbox{}\hspace{8pt}\tau_{52}=0\mbox{}\hspace{8pt}\tau_{61}=1\mbox{}\hspace{8pt}\beta_{6}=\tau_{62}=0\\
\beta_{7}=\tau_{71}=1\mbox{}\hspace{8pt}\tau_{72}=0\mbox{}\hspace{8pt}\tau_{81}=1\mbox{}\hspace{8pt}\beta_{8}=\tau_{82}=0\\
\beta_{9}=\tau_{91}=1\mbox{}\hspace{8pt}\tau_{92}=0\mbox{}\hspace{8pt}\tau_{(10)1}=1\mbox{}\hspace{8pt}\beta_{10}=\tau_{(10)2}=0\\
\beta_{11}=\tau_{(11)1}=1\mbox{}\hspace{8pt}\tau_{(11)2}=0\mbox{}\hspace{8pt}\tau_{(12)1}=1\mbox{}\hspace{8pt}\beta_{12}=\tau_{(12)2}=0\\
\beta_{13}=\tau_{(13)1}=1\mbox{}\hspace{8pt}\tau_{(13)2}=0\\
\end{array}
\right.
$$
Substituting $\beta_{13}=\tau_{(13)1}=1$ and $\tau_{(12)2}=0$
into the fourteen equation of Eqs.(48) yields
$$3=\tau_{(14)1}2+\tau_{(14)2}2^{2}$$
which has no solution.
\par It is not hard to verify $g(\gamma_{14})=1$, by Lemma 7, Eqs.(48) has no solution.
\par {\bf Case 2.} $(m,n)=(2,13)$, we have
$$2^{16}+1=(2^{2}+\alpha_{1}2+1)(2^{13}+\sum\limits^{12}_{i=1}\beta_{i}2^{i}+1)$$
Expanding the right hand side of the above equation yields
$$
\begin{array}{lll}
2^{16}+1&=&2^{15}\\
&+&(\alpha_{1}+\beta_{12})2^{14}\\
&+&(1+\alpha_{1}\beta_{12}+\beta_{11})2^{13}\\
&+&\sum\limits^{12}_{i=3}(\beta_{i}+\alpha_{1}\beta_{i-1}+\beta_{i-2})2^{i}\\
&+&(\beta_{2}+\alpha_{1}\beta_{1}+1)2^{2}\\
&+&(\beta_{1}+\alpha_{1})2\\
&+&1
\end{array}
$$
which implies that
\begin{equation}
\left\{
\begin{array}{lr}
\beta_{1}+\alpha_{1}=\tau_{11}2\\
\beta_{2}+\alpha_{1}\beta_{1}+1+\tau_{11}=\tau_{21}2+\tau_{22}2^{2}\\
\beta_{3}+\alpha_{1}\beta_{2}+\beta_{1}+\tau_{21}=\tau_{31}2+\tau_{32}2^{2}\\
\beta_{i}+\alpha_{1}\beta_{i-1}+\beta_{i-2}+\tau_{(i-1)1}+\tau_{(i-2)2}=\tau_{i1}2+\tau_{i2}2^{2}
\mbox{}\hspace{12pt} 4\leq i\leq 12\\
1+\alpha_{1}\beta_{12}+\beta_{11}+\tau_{(12)1}+\tau_{(11)2}=\tau_{(13)1}2+\tau_{(13)2}2^{2}\\
\alpha_{1}+\beta_{12}+\tau_{(13)1}+\tau_{(12)2}=\tau_{(14)1}2+\tau_{(14)2}2^{2}\\
1+\tau_{(14)1}+\tau_{(13)2}=\tau_{(15)1}2\\
\tau_{(15)1}+\tau_{(14)2}=1
\end{array}
\right.
\end{equation}
\par There are two situations.
\par {\bf Case 2-1.} If $\alpha_{1}=0$,
then Eqs.(49) become
\begin{equation}
\left\{
\begin{array}{lr}
\beta_{1}=\tau_{11}2\\
\beta_{2}+1+\tau_{11}=\tau_{21}2+\tau_{22}2^{2}\\
\beta_{3}+\beta_{1}+\tau_{21}=\tau_{31}2+\tau_{32}2^{2}\\
\beta_{i}+\beta_{i-2}+\tau_{(i-1)1}+\tau_{(i-2)2}=\tau_{i1}2+\tau_{i2}2^{2}
\mbox{}\hspace{12pt} 4\leq i\leq 12\\
1+\alpha_{1}\beta_{12}+\beta_{11}+\tau_{(12)1}+\tau_{(11)2}=\tau_{(13)1}2+\tau_{(13)2}2^{2}\\
\alpha_{1}+\beta_{12}+\tau_{(13)1}+\tau_{(12)2}=\tau_{(14)1}2+\tau_{(14)2}2^{2}\\
1+\tau_{(14)1}+\tau_{(13)2}=\tau_{(15)1}2\\
\tau_{(15)1}+\tau_{(14)2}=1
\end{array}
\right.
\end{equation}
By Lemma 6, we have $g(\gamma_{14})=1$. Thus, Eqs.(50) has no solution.
\par {\bf Case 2-2.} If $\alpha_{1}=1$,
then Eqs.(49) become
\begin{equation}
\left\{
\begin{array}{lr}
\beta_{1}+1=\tau_{11}2\\
\beta_{2}+\beta_{1}+1+\tau_{11}=\tau_{21}2+\tau_{22}2^{2}\\
\beta_{3}+\beta_{2}+\beta_{1}+\tau_{21}=\tau_{31}2+\tau_{32}2^{2}\\
\beta_{i}+\beta_{i-1}+\beta_{i-2}+\tau_{(i-1)1}+\tau_{(i-2)2}=\tau_{i1}2+\tau_{i2}2^{2}
\mbox{}\hspace{12pt} 4\leq i\leq 12\\
1+\beta_{12}+\beta_{11}+\tau_{(12)1}+\tau_{(11)2}=\tau_{(13)1}2+\tau_{(13)2}2^{2}\\
\alpha_{1}+\beta_{12}+\tau_{(13)1}+\tau_{(12)2}=\tau_{(14)1}2+\tau_{(14)2}2^{2}\\
1+\tau_{(14)1}+\tau_{(13)2}=\tau_{(15)1}2\\
\tau_{(15)1}+\tau_{(14)2}=1
\end{array}
\right.
\end{equation}
Similarly to case 1, we have $g(\gamma_{14})=1$. Thus, Eqs.(51) has no solution.
\par {\bf Case 3.} $(m,n)=(3,12)$, we have
$$2^{16}+1=(2^{3}+\alpha_{2}2^{2}+\alpha_{1}2+1)(2^{12}+\sum\limits^{11}_{i=1}\beta_{i}2^{i}+1)$$
The equations satisfied by variables $\alpha_{i}(i=1,2)$ and
$\beta_{j}(j=1,2,\cdots,11)$ are given by
$$
\left\{
\begin{array}{lr}
\beta_{1}+\alpha_{1}=\tau_{11}2\\
\alpha_{2}+\alpha_{1}\beta_{1}+\beta_{2}+\tau_{11}=\tau_{21}2+\tau_{22}2^{2}\\
\beta_{3}+\alpha_{1}\beta_{2}+\alpha_{2}\beta_{1}+1+\tau_{21}=\tau_{31}2+\tau_{32}2^{2}\\
\beta_{i}+\alpha_{1}\beta_{i-1}+\alpha_{2}\beta_{i-2}+\beta_{i-3}+\tau_{(i-1)1}+\tau_{(i-2)2}=\tau_{i1}2+\tau_{i2}2^{2}
\mbox{}\hspace{8pt} 4\leq i\leq 11\\
1+\alpha_{1}\beta_{11}+\alpha_{2}\beta_{10}+\beta_{9}+\tau_{(11)1}+\tau_{(10)2}=\tau_{(12)1}2+\tau_{(12)2}2^{2}\\
\alpha_{1}+\alpha_{2}\beta_{11}+\beta_{10}+\tau_{(12)1}+\tau_{(11)2}=\tau_{(13)1}2+\tau_{(13)2}2^{2}\\
\alpha_{2}+\beta_{11}+\tau_{(13)1}+\tau_{(12)2}=\tau_{(14)1}2+\tau_{(14)2}2^{2}\\
1+\tau_{(14)1}+\tau_{(13)2}=\tau_{(15)1}2\\
\tau_{(15)1}+\tau_{(14)2}=1
\end{array}
\right.
$$
\par There are four cases.
\par {\bf Case 3-1.} $\alpha_{1}=\alpha_{2}=0$. In this
case, we get
$$
\left\{
\begin{array}{lr}
\beta_{1}=\tau_{11}2\\
\alpha_{2}+\beta_{2}+\tau_{11}=\tau_{21}2+\tau_{22}2^{2}\\
\beta_{3}+1+\tau_{21}=\tau_{31}2+\tau_{32}2^{2}\\
\beta_{i}+\beta_{i-3}+\tau_{(i-1)1}+\tau_{(i-2)2}=\tau_{i1}2+\tau_{i2}2^{2}
\mbox{}\hspace{8pt} 4\leq i\leq 11\\
1+\beta_{9}+\tau_{(11)1}+\tau_{(10)2}=\tau_{(12)1}2+\tau_{(12)2}2^{2}\\
\beta_{10}+\tau_{(12)1}+\tau_{(11)2}=\tau_{(13)1}2+\tau_{(13)2}2^{2}\\
\beta_{11}+\tau_{(13)1}+\tau_{(12)2}=\tau_{(14)1}2+\tau_{(14)2}2^{2}\\
1+\tau_{(14)1}+\tau_{(13)2}=\tau_{(15)1}2\\
\tau_{(15)1}+\tau_{(14)2}=1
\end{array}
\right.
$$
By Lemma 6, we have $g(\gamma_{12})=1$. Thus, the above equations has no solution.
\par {\bf Case 3-2.}  $\alpha_{1}=0$ and $\alpha_{2}=1$. In
this case, we have
$$
\left\{
\begin{array}{lr}
\beta_{1}=\tau_{11}2\\
\alpha_{2}+\beta_{2}+\tau_{11}=\tau_{21}2+\tau_{22}2^{2}\\
\beta_{3}+\beta_{1}+1+\tau_{21}=\tau_{31}2+\tau_{32}2^{2}\\
\beta_{i}+\beta_{i-2}+\beta_{i-3}+\tau_{(i-1)1}+\tau_{(i-2)2}=\tau_{i1}2+\tau_{i2}2^{2}
\mbox{}\hspace{8pt} 4\leq i\leq 11\\
1+\beta_{10}+\beta_{9}+\tau_{(11)1}+\tau_{(10)2}=\tau_{(12)1}2+\tau_{(12)2}2^{2}\\
\beta_{11}+\beta_{10}+\tau_{(12)1}+\tau_{(11)2}=\tau_{(13)1}2+\tau_{(13)2}2^{2}\\
1+\beta_{11}+\tau_{(13)1}+\tau_{(12)2}=\tau_{(14)1}2+\tau_{(14)2}2^{2}\\
1+\tau_{(14)1}+\tau_{(13)2}=\tau_{(15)1}2\\
\tau_{(15)1}+\tau_{(14)2}=1
\end{array}
\right.
$$
By Lemma 6, we have $g(\gamma_{12})=1$. Thus, the above equations has no solution.
\par {\bf Case 3-3.}
$\alpha_{1}=1$ and $\alpha_{2}=0$. In
this case, we obtain
$$
\left\{
\begin{array}{lr}
\beta_{1}+1=\tau_{11}2\\
\beta_{1}+\beta_{2}+\tau_{11}=\tau_{21}2+\tau_{22}2^{2}\\
\beta_{3}+\beta_{2}+1+\tau_{21}=\tau_{31}2+\tau_{32}2^{2}\\
\beta_{i}+\beta_{i-1}+\beta_{i-3}+\tau_{(i-1)1}+\tau_{(i-2)2}=\tau_{i1}2+\tau_{i2}2^{2}
\mbox{}\hspace{8pt} 4\leq i\leq 11\\
1+\beta_{11}+\beta_{9}+\tau_{(11)1}+\tau_{(10)2}=\tau_{(12)1}2+\tau_{(12)2}2^{2}\\
1+\beta_{10}+\tau_{(12)1}+\tau_{(11)2}=\tau_{(13)1}2+\tau_{(13)2}2^{2}\\
\beta_{11}+\tau_{(13)1}+\tau_{(12)2}=\tau_{(14)1}2+\tau_{(14)2}2^{2}\\
1+\tau_{(14)1}+\tau_{(13)2}=\tau_{(15)1}2\\
\tau_{(15)1}+\tau_{(14)2}=1
\end{array}
\right.
$$
By Lemma 6, we have $g(\gamma_{12})=1$. Thus, the above equations has no solution.
\par {\bf Case 3-4.} $\alpha_{1}=\alpha_{2}=1$. In
this case, we get
$$
\left\{
\begin{array}{lr}
\beta_{1}+1=\tau_{11}2\\
1+\beta_{1}+\beta_{2}+\tau_{11}=\tau_{21}2+\tau_{22}2^{2}\\
\beta_{3}+\beta_{2}+\beta_{1}+1+\tau_{21}=\tau_{31}2+\tau_{32}2^{2}\\
\beta_{i}+\beta_{i-1}+\beta_{i-2}+\beta_{i-3}+\tau_{(i-1)1}+\tau_{(i-2)2}=\tau_{i1}2+\tau_{i2}2^{2}
\mbox{}\hspace{8pt} 4\leq i\leq 11\\
1+\beta_{11}+\beta_{10}+\beta_{9}+\tau_{(11)1}+\tau_{(10)2}=\tau_{(12)1}2+\tau_{(12)2}2^{2}\\
1+\beta_{11}+\beta_{10}+\tau_{(12)1}+\tau_{(11)2}=\tau_{(13)1}2+\tau_{(13)2}2^{2}\\
1+\beta_{11}+\tau_{(13)1}+\tau_{(12)2}=\tau_{(14)1}2+\tau_{(14)2}2^{2}\\
1+\tau_{(14)1}+\tau_{(13)2}=\tau_{(15)1}2\\
\tau_{(15)1}+\tau_{(14)2}=1
\end{array}
\right.
$$
By Lemma 6, we have $g(\gamma_{12})=1$. Thus, the above equations has no solution.
\par {\bf Case 4.} $(m,n)=(4,11)$, the equations satisfied by variables $\alpha_{i}(i=1,2,3)$ and
$\beta_{j}(j=1,2,\cdots,10)$ are given by
$$
\left\{
\begin{array}{lr}
\beta_{1}+\alpha_{1}=\tau_{11}2\\
\beta_{2}+\alpha_{1}\beta_{1}+\alpha_{2}+\tau_{11}=\tau_{21}2+\tau_{22}2^{2}\\
\beta_{3}+\alpha_{1}\beta_{2}+\alpha_{2}\beta_{1}+\alpha_{3}+\tau_{21}=\tau_{31}2+\tau_{32}2^{2}\\
\beta_{4}+\alpha_{1}\beta_{3}+\alpha_{2}\beta_{2}+\alpha_{3}\beta_{1}+1+\tau_{31}+\tau_{22}=\tau_{41}2+\tau_{42}2^{2}\\
\beta_{i}+\alpha_{1}\beta_{i-1}+\alpha_{2}\beta_{i-2}+\alpha_{3}\beta_{i-3}+\beta_{i-4}+\tau_{(i-1)1}+\tau_{(i-2)2}=\tau_{i1}2+\tau_{i2}2^{2}\\
1+\alpha_{1}\beta_{10}+\alpha_{2}\beta_{9}+\alpha_{3}\beta_{8}+\beta_{7}+\tau_{(10)1}+\tau_{92}=\tau_{(11)1}2+\tau_{(11)2}2^{2}\\
\alpha_{1}+\alpha_{2}\beta_{10}+\alpha_{3}\beta_{9}+\beta_{8}+\tau_{(11)1}+\tau_{(10)2}=\tau_{(12)1}2+\tau_{(12)2}2^{2}\\
\alpha_{2}+\alpha_{3}\beta_{10}+\beta_{9}+\tau_{(12)1}+\tau_{(11)2}=\tau_{(13)1}2+\tau_{(13)2}2^{2}\\
\alpha_{3}+\beta_{10}+\tau_{(13)1}+\tau_{(12)2}=\tau_{(14)1}2+\tau_{(14)2}2^{2}\\
1+\tau_{(14)1}+\tau_{(13)2}=\tau_{(15)1}2\\
\tau_{(15)1}+\tau_{(14)2}=1
\end{array}
\right.
$$
where $5\leq i\leq 10$.
\par {\bf Case 4-1.} $\alpha_{1}=\alpha_{2}=\alpha_{3}=0$. We get
$$
\left\{
\begin{array}{lr}
\beta_{1}=\tau_{11}2\\
\beta_{2}+\tau_{11}=\tau_{21}2+\tau_{22}2^{2}\\
\beta_{3}+\alpha_{1}\beta_{2}+\alpha_{2}\beta_{1}+\alpha_{3}+\tau_{21}=\tau_{31}2+\tau_{32}2^{2}\\
\beta_{4}+1+\tau_{31}+\tau_{22}=\tau_{41}2+\tau_{42}2^{2}\\
\beta_{i}+\beta_{i-4}+\tau_{(i-1)1}+\tau_{(i-2)2}=\tau_{i1}2+\tau_{i2}2^{2}\mbox{}\hspace{12pt} 5\leq i\leq 10\\
1+\beta_{7}+\tau_{(10)1}+\tau_{92}=\tau_{(11)1}2+\tau_{(11)2}2^{2}\\
\beta_{8}+\tau_{(11)1}+\tau_{(10)2}=\tau_{(12)1}2+\tau_{(12)2}2^{2}\\
\beta_{9}+\tau_{(12)1}+\tau_{(11)2}=\tau_{(13)1}2+\tau_{(13)2}2^{2}\\
\beta_{10}+\tau_{(13)1}+\tau_{(12)2}=\tau_{(14)1}2+\tau_{(14)2}2^{2}\\
1+\tau_{(14)1}+\tau_{(13)2}=\tau_{(15)1}2\\
\tau_{(15)1}+\tau_{(14)2}=1
\end{array}
\right.
$$
In this case, $g(\gamma_{11})=1$.
\par {\bf Case 4-2.} $\alpha_{1}=\alpha_{2}=0$ and $\alpha_{3}=1$. We obtain
$$
\left\{
\begin{array}{lr}
\beta_{1}=\tau_{11}2\\
\beta_{2}+\tau_{11}=\tau_{21}2+\tau_{22}2^{2}\\
\beta_{3}+1+\tau_{21}=\tau_{31}2+\tau_{32}2^{2}\\
\beta_{4}+\beta_{1}+1+\tau_{31}+\tau_{22}=\tau_{41}2+\tau_{42}2^{2}\\
\beta_{i}+\beta_{i-3}+\beta_{i-4}+\tau_{(i-1)1}+\tau_{(i-2)2}=\tau_{i1}2+\tau_{i2}2^{2}\mbox{}\hspace{12pt} 5\leq i\leq 10\\
1+\beta_{8}+\beta_{7}+\tau_{(10)1}+\tau_{92}=\tau_{(11)1}2+\tau_{(11)2}2^{2}\\
\beta_{9}+\beta_{8}+\tau_{(11)1}+\tau_{(10)2}=\tau_{(12)1}2+\tau_{(12)2}2^{2}\\
\beta_{10}+\beta_{9}+\tau_{(12)1}+\tau_{(11)2}=\tau_{(13)1}2+\tau_{(13)2}2^{2}\\
1+\beta_{10}+\tau_{(13)1}+\tau_{(12)2}=\tau_{(14)1}2+\tau_{(14)2}2^{2}\\
1+\tau_{(14)1}+\tau_{(13)2}=\tau_{(15)1}2\\
\tau_{(15)1}+\tau_{(14)2}=1
\end{array}
\right.
$$
In this case, $g(\gamma_{12})=1$.
\par {\bf Case 4-3.} $\alpha_{1}=0$, $\alpha_{2}=1$ and $\alpha_{3}=0$. We have
$$
\left\{
\begin{array}{lr}
\beta_{1}=\tau_{11}2\\
\beta_{2}+1+\tau_{11}=\tau_{21}2+\tau_{22}2^{2}\\
\beta_{3}+\beta_{1}+\tau_{21}=\tau_{31}2+\tau_{32}2^{2}\\
\beta_{4}+\beta_{2}+1+\tau_{31}+\tau_{22}=\tau_{41}2+\tau_{42}2^{2}\\
\beta_{i}+\beta_{i-2}+\beta_{i-4}+\tau_{(i-1)1}+\tau_{(i-2)2}=\tau_{i1}2+\tau_{i2}2^{2}\mbox{}\hspace{12pt} 5\leq i\leq 10\\
1+\beta_{9}+\beta_{7}+\tau_{(10)1}+\tau_{92}=\tau_{(11)1}2+\tau_{(11)2}2^{2}\\
\beta_{10}+\beta_{8}+\tau_{(11)1}+\tau_{(10)2}=\tau_{(12)1}2+\tau_{(12)2}2^{2}\\
1+\beta_{9}+\tau_{(12)1}+\tau_{(11)2}=\tau_{(13)1}2+\tau_{(13)2}2^{2}\\
\beta_{10}+\tau_{(13)1}+\tau_{(12)2}=\tau_{(14)1}2+\tau_{(14)2}2^{2}\\
1+\tau_{(14)1}+\tau_{(13)2}=\tau_{(15)1}2\\
\tau_{(15)1}+\tau_{(14)2}=1
\end{array}
\right.
$$
In this case, $g(\gamma_{14})=1$.
\par {\bf Case 4-4.} $\alpha_{1}=1$ and $\alpha_{2}=\alpha_{3}=0$. We get
$$
\left\{
\begin{array}{lr}
\beta_{1}+1=\tau_{11}2\\
\beta_{2}+\beta_{1}+\tau_{11}=\tau_{21}2+\tau_{22}2^{2}\\
\beta_{3}+\beta_{2}+\tau_{21}=\tau_{31}2+\tau_{32}2^{2}\\
\beta_{4}+\beta_{3}+1+\tau_{31}+\tau_{22}=\tau_{41}2+\tau_{42}2^{2}\\
\beta_{i}+\beta_{i-1}+\beta_{i-4}+\tau_{(i-1)1}+\tau_{(i-2)2}=\tau_{i1}2+\tau_{i2}2^{2}\mbox{}\hspace{12pt} 5\leq i\leq 10\\
1+\beta_{10}+\beta_{7}+\tau_{(10)1}+\tau_{92}=\tau_{(11)1}2+\tau_{(11)2}2^{2}\\
1+\beta_{8}+\tau_{(11)1}+\tau_{(10)2}=\tau_{(12)1}2+\tau_{(12)2}2^{2}\\
\beta_{9}+\tau_{(12)1}+\tau_{(11)2}=\tau_{(13)1}2+\tau_{(13)2}2^{2}\\
\beta_{10}+\tau_{(13)1}+\tau_{(12)2}=\tau_{(14)1}2+\tau_{(14)2}2^{2}\\
1+\tau_{(14)1}+\tau_{(13)2}=\tau_{(15)1}2\\
\tau_{(15)1}+\tau_{(14)2}=1
\end{array}
\right.
$$
In this case, $g(\gamma_{14})=1$.
\par {\bf Case 4-5.} $\alpha_{1}=0$ and $\alpha_{2}=\alpha_{3}=1$. We have
$$
\left\{
\begin{array}{lr}
\beta_{1}=\tau_{11}2\\
\beta_{2}+1+\tau_{11}=\tau_{21}2+\tau_{22}2^{2}\\
\beta_{3}+\beta_{1}+1+\tau_{21}=\tau_{31}2+\tau_{32}2^{2}\\
\beta_{4}+\beta_{2}+\beta_{1}+1+\tau_{31}+\tau_{22}=\tau_{41}2+\tau_{42}2^{2}\\
\beta_{i}+\beta_{i-2}+\beta_{i-3}+\beta_{i-4}+\tau_{(i-1)1}+\tau_{(i-2)2}=\tau_{i1}2+\tau_{i2}2^{2}\mbox{}\hspace{12pt} 5\leq i\leq 10\\
1+\beta_{9}+\beta_{8}+\beta_{7}+\tau_{(10)1}+\tau_{92}=\tau_{(11)1}2+\tau_{(11)2}2^{2}\\
\beta_{10}+\beta_{9}+\beta_{8}+\tau_{(11)1}+\tau_{(10)2}=\tau_{(12)1}2+\tau_{(12)2}2^{2}\\
1+\beta_{10}+\beta_{9}+\tau_{(12)1}+\tau_{(11)2}=\tau_{(13)1}2+\tau_{(13)2}2^{2}\\
1+\beta_{10}+\tau_{(13)1}+\tau_{(12)2}=\tau_{(14)1}2+\tau_{(14)2}2^{2}\\
1+\tau_{(14)1}+\tau_{(13)2}=\tau_{(15)1}2\\
\tau_{(15)1}+\tau_{(14)2}=1
\end{array}
\right.
$$
In this case, $g(\gamma_{11})=1$.
\par {\bf Case 4-6.} $\alpha_{1}=1$, $\alpha_{2}=0$ and $\alpha_{3}=1$.
$$
\left\{
\begin{array}{lr}
\beta_{1}+1=\tau_{11}2\\
\beta_{2}+\beta_{1}+\tau_{11}=\tau_{21}2+\tau_{22}2^{2}\\
\beta_{3}+\beta_{2}+1+\tau_{21}=\tau_{31}2+\tau_{32}2^{2}\\
\beta_{4}+\beta_{3}+\beta_{1}+1+\tau_{31}+\tau_{22}=\tau_{41}2+\tau_{42}2^{2}\\
\beta_{i}+\beta_{i-1}+\beta_{i-3}+\beta_{i-4}+\tau_{(i-1)1}+\tau_{(i-2)2}=\tau_{i1}2+\tau_{i2}2^{2}\mbox{}\hspace{12pt} 5\leq i\leq 10\\
1+\beta_{10}+\beta_{8}+\beta_{7}+\tau_{(10)1}+\tau_{92}=\tau_{(11)1}2+\tau_{(11)2}2^{2}\\
1+\beta_{9}+\beta_{8}+\tau_{(11)1}+\tau_{(10)2}=\tau_{(12)1}2+\tau_{(12)2}2^{2}\\
\beta_{10}+\beta_{9}+\tau_{(12)1}+\tau_{(11)2}=\tau_{(13)1}2+\tau_{(13)2}2^{2}\\
1+\beta_{10}+\tau_{(13)1}+\tau_{(12)2}=\tau_{(14)1}2+\tau_{(14)2}2^{2}\\
1+\tau_{(14)1}+\tau_{(13)2}=\tau_{(15)1}2\\
\tau_{(15)1}+\tau_{(14)2}=1
\end{array}
\right.
$$
In this case, $g(\gamma_{12})=1$.
\par {\bf Case 4-7.} $\alpha_{1}=\alpha_{2}=1$ and $\alpha_{3}=0$. We have
$$
\left\{
\begin{array}{lr}
\beta_{1}+1=\tau_{11}2\\
\beta_{2}+\beta_{1}+1+\tau_{11}=\tau_{21}2+\tau_{22}2^{2}\\
\beta_{3}+\beta_{2}+\beta_{1}+\tau_{21}=\tau_{31}2+\tau_{32}2^{2}\\
\beta_{4}+\beta_{3}+\beta_{2}+1+\tau_{31}+\tau_{22}=\tau_{41}2+\tau_{42}2^{2}\\
\beta_{i}+\beta_{i-1}+\beta_{i-2}+\beta_{i-4}+\tau_{(i-1)1}+\tau_{(i-2)2}=\tau_{i1}2+\tau_{i2}2^{2}\mbox{}\hspace{12pt} 5\leq i\leq 10\\
1+\beta_{10}+\beta_{9}+\beta_{7}+\tau_{(10)1}+\tau_{92}=\tau_{(11)1}2+\tau_{(11)2}2^{2}\\
1+\beta_{10}+\beta_{8}+\tau_{(11)1}+\tau_{(10)2}=\tau_{(12)1}2+\tau_{(12)2}2^{2}\\
1+\beta_{9}+\tau_{(12)1}+\tau_{(11)2}=\tau_{(13)1}2+\tau_{(13)2}2^{2}\\
\beta_{10}+\tau_{(13)1}+\tau_{(12)2}=\tau_{(14)1}2+\tau_{(14)2}2^{2}\\
1+\tau_{(14)1}+\tau_{(13)2}=\tau_{(15)1}2\\
\tau_{(15)1}+\tau_{(14)2}=1
\end{array}
\right.
$$
In this case, $g(\gamma_{11})=1$.
\par {\bf Case 4-8.} $\alpha_{1}=\alpha_{2}=\alpha_{3}=1$. We find that
$$
\left\{
\begin{array}{lr}
\beta_{1}+1=\tau_{11}2\\
\beta_{2}+\beta_{1}+1+\tau_{11}=\tau_{21}2+\tau_{22}2^{2}\\
\beta_{3}+\beta_{2}+\beta_{1}+1+\tau_{21}=\tau_{31}2+\tau_{32}2^{2}\\
\beta_{4}+\beta_{3}+\beta_{2}+\beta_{1}+1+\tau_{31}+\tau_{22}=\tau_{41}2+\tau_{42}2^{2}\\
\beta_{i}+\beta_{i-1}+\beta_{i-2}+\beta_{i-3}+\beta_{i-4}+\tau_{(i-1)1}+\tau_{(i-2)2}=\tau_{i1}2+\tau_{i2}2^{2}\\
1+\beta_{10}+\beta_{9}+\beta_{8}+\beta_{7}+\tau_{(10)1}+\tau_{92}=\tau_{(11)1}2+\tau_{(11)2}2^{2}\\
1+\beta_{10}+\beta_{9}+\beta_{8}+\tau_{(11)1}+\tau_{(10)2}=\tau_{(12)1}2+\tau_{(12)2}2^{2}\\
1+\beta_{10}+\beta_{9}+\tau_{(12)1}+\tau_{(11)2}=\tau_{(13)1}2+\tau_{(13)2}2^{2}\\
1+\beta_{10}+\tau_{(13)1}+\tau_{(12)2}=\tau_{(14)1}2+\tau_{(14)2}2^{2}\\
1+\tau_{(14)1}+\tau_{(13)2}=\tau_{(15)1}2\\
\tau_{(15)1}+\tau_{(14)2}=1
\end{array}
\right.
$$
In this case, $g(\gamma_{12})=1$.

\par {\bf Case 5.} $(m,n)=(5,10)$, the equations satisfied by variables $\alpha_{i}(i=1,2,3,4)$ and
$\beta_{j}(j=1,2,\cdots,9)$ are given by
$$
\left\{
\begin{array}{lr}
\beta_{1}+\alpha_{1}=\tau_{11}2\\
\beta_{2}+\alpha_{1}\beta_{1}+\alpha_{2}+\tau_{11}=\tau_{21}2+\tau_{22}2^{2}\\
\beta_{3}+\alpha_{1}\beta_{2}+\alpha_{2}\beta_{1}+\alpha_{3}+\tau_{21}=\tau_{31}2+\tau_{32}2^{2}\\
\beta_{4}+\alpha_{1}\beta_{3}+\alpha_{2}\beta_{2}+\alpha_{3}\beta_{1}+\alpha_{4}+\tau_{31}+\tau_{22}=\tau_{41}2+\tau_{42}2^{2}\\
\beta_{5}+\alpha_{1}\beta_{4}+\alpha_{2}\beta_{3}+\alpha_{3}\beta_{2}+\alpha_{4}\beta_{1}+1+\tau_{41}+\tau_{32}=\tau_{51}2+\tau_{52}2^{2}+\tau_{53}2^{3}\\
\beta_{6}+\alpha_{1}\beta_{5}+\alpha_{2}\beta_{4}+\alpha_{3}\beta_{3}+\alpha_{4}\beta_{2}+\beta_{1}+\tau_{51}++\tau_{42}=\tau_{61}2+\tau_{62}2^{2}+\tau_{63}2^{3}\\
\beta_{7}+\alpha_{1}\beta_{6}+\alpha_{2}\beta_{5}+\alpha_{3}\beta_{4}+\alpha_{4}\beta_{3}+\beta_{2}+\tau_{61}+\tau_{52}=\tau_{71}2+\tau_{72}2^{2}+\tau_{73}2^{3}\\
\beta_{8}+\sum\limits^{4}_{i=1}\alpha_{i}\beta_{8-i}+\beta_{3}+\tau_{71}+\tau_{62}+\tau_{53}=\tau_{81}2+\tau_{82}2^{2}+\tau_{83}2^{3}\\
\beta_{9}+\sum\limits^{4}_{i=1}\alpha_{i}\beta_{9-i}+\beta_{4}+\tau_{81}+\tau_{72}+\tau_{63}=\tau_{91}2+\tau_{92}2^{2}+\tau_{93}2^{3}\\
1+\sum\limits^{4}_{i=1}\alpha_{i}\beta_{10-i}+\beta_{5}+\tau_{91}+\tau_{82}+\tau_{73}=\tau_{(10)1}2+\tau_{(10)2}2^{2}+\tau_{(10)3}2^{3}\\
\alpha_{1}+\alpha_{2}\beta_{9}+\alpha_{3}\beta_{8}+\alpha_{4}\beta_{7}+\beta_{6}+\tau_{(10)1}+\tau_{92}+\tau_{83}=\tau_{(11)1}2+\tau_{(11)2}2^{2}+\tau_{(11)3}2^{3}\\
\alpha_{2}+\alpha_{3}\beta_{9}+\alpha_{4}\beta_{8}+\beta_{7}+\tau_{(11)1}+\tau_{(10)2}+\tau_{93}=\tau_{(12)1}2+\tau_{(12)2}2^{2}\\
\alpha_{3}+\alpha_{4}\beta_{9}+\beta_{8}+\tau_{(12)1}+\tau_{(11)2}+\tau_{(10)3}=\tau_{(13)1}2+\tau_{(13)2}2^{2}\\
\alpha_{4}+\beta_{9}+\tau_{(13)1}+\tau_{(12)2}+\tau_{(11)3}=\tau_{(14)1}2+\tau_{(14)2}2^{2}\\
1+\tau_{(14)1}+\tau_{(13)2}=\tau_{(15)1}2\\
\tau_{(15)1}+\tau_{(14)2}=1
\end{array}
\right.
$$
\par There are sixteen cases.
\par {\bf Case 5-1.} $\alpha_{1}=\alpha_{2}=\alpha_{3}=\alpha_{4}=0$. In this case, $g(\gamma_{10})=1$.
\par {\bf Case 5-2.} $\alpha_{1}=\alpha_{2}=\alpha_{3}=0$ and $\alpha_{4}=1$. In this case, $g(\gamma_{10})=1$.
\par {\bf Case 5-3.} $\alpha_{1}=\alpha_{2}=0$, $\alpha_{3}=1$ and $\alpha_{4}=0$. In this case, $g(\gamma_{11})=1$.
\par {\bf Case 5-4.} $\alpha_{1}=\alpha_{2}=0$ and $\alpha_{3}=\alpha_{4}=0$. In this case, $g(\gamma_{11})=1$.
\par {\bf Case 5-5.} $\alpha_{1}=0$, $\alpha_{2}=1$ and $\alpha_{3}=\alpha_{4}=0$. In this case, $g(\gamma_{10})=1$.
\par {\bf Case 5-6.} $\alpha_{1}=0$, $\alpha_{2}=1$, $\alpha_{3}=0$ and $\alpha_{4}=1$. In this case, $g(\gamma_{10})=1$.
\par {\bf Case 5-7.} $\alpha_{1}=0$, $\alpha_{2}=\alpha_{3}=1$ and $\alpha_{4}=0$. In this case, $g(\gamma_{11})=1$.
\par {\bf Case 5-8.} $\alpha_{1}=0$ and $\alpha_{2}=\alpha_{3}=\alpha_{4}=1$. In this case, $g(\gamma_{12})=1$.
\par {\bf Case 5-9.} $\alpha_{1}=1$ and $\alpha_{2}=\alpha_{3}=\alpha_{4}=0$. In this case, $g(\gamma_{10})=1$.
\par {\bf Case 5-10.} $\alpha_{1}=1$, $\alpha_{2}=\alpha_{3}=0$ and $\alpha_{4}=0$. In this case, $g(\gamma_{10})=1$.
\par {\bf Case 5-11.} $\alpha_{1}=1$, $\alpha_{2}=0$, $\alpha_{3}=1$ and $\alpha_{4}=0$. In this case, $g(\gamma_{11})=1$.
\par {\bf Case 5-12.} $\alpha_{1}=1$, $\alpha_{2}=0$ and $\alpha_{3}=\alpha_{4}=1$. In this case, $g(\gamma_{10})=1$.
\par {\bf Case 5-13.} $\alpha_{1}=\alpha_{2}=1$ and $\alpha_{3}=\alpha_{4}=0$. In this case, $g(\gamma_{10})=1$.
\par {\bf Case 5-14.} $\alpha_{1}=\alpha_{2}=1$, $\alpha_{3}=0$ and $\alpha_{4}=1$. In this case, $g(\gamma_{11})=1$.
\par {\bf Case 5-15.} $\alpha_{1}=\alpha_{2}=\alpha_{3}=1$ and $\alpha_{4}=0$. In this case, $g(\gamma_{14})=1$.
\par {\bf Case 5-16.} $\alpha_{1}=\alpha_{2}=\alpha_{3}=\alpha_{4}=1$. In this case, $g(\gamma_{11})=1$.
\par {\bf Case 6.} $(m,n)=(6,9)$. There are thirty-two cases.
\par {\bf Case 6-1.} $\alpha_{1}=\alpha_{2}=\alpha_{3}=\alpha_{4}=\alpha_{5}=0$. In this case, $g(\gamma_{10})=1$.
\par {\bf Case 6-2.} $\alpha_{1}=\alpha_{2}=\alpha_{3}=\alpha_{4}=0$ and $\alpha_{5}=1$. In this case, $g(\gamma_{13})=1$.
\par {\bf Case 6-3.} $\alpha_{1}=\alpha_{2}=\alpha_{3}=0$, $\alpha_{4}=1$ and $\alpha_{5}=0$. In this case, $g(\gamma_{9})=1$.
\par {\bf Case 6-4.} $\alpha_{1}=\alpha_{2}=\alpha_{3}=0$ and $\alpha_{4}=\alpha_{5}=1$. In this case, $g(\gamma_{9})=1$.
\par {\bf Case 6-5.} $\alpha_{1}=\alpha_{2}=0$, $\alpha_{3}=1$ and $\alpha_{4}=\alpha_{5}=0$. In this case, $g(\gamma_{9})=1$.
\par {\bf Case 6-6.} $\alpha_{1}=\alpha_{2}=0$, $\alpha_{3}=1$, $\alpha_{4}=0$ and $\alpha_{5}=1$. In this case, $g(\gamma_{11})=1$.
\par {\bf Case 6-7.} $\alpha_{1}=\alpha_{2}=0$, $\alpha_{3}=\alpha_{4}=1$ and $\alpha_{5}=0$. In this case, $g(\gamma_{10})=1$.
\par {\bf Case 6-8.} $\alpha_{1}=\alpha_{2}=0$ and $\alpha_{3}=\alpha_{4}=\alpha_{5}=1$. In this case, $g(\gamma_{9})=1$.
\par {\bf Case 6-9.} $\alpha_{1}=0$, $\alpha_{2}=1$ and $\alpha_{3}=\alpha_{4}=\alpha_{5}=0$. In this case, $g(\gamma_{10})=1$.
\par {\bf Case 6-10.} $\alpha_{1}=0$, $\alpha_{2}=1$, $\alpha_{3}=\alpha_{4}=0$ and $\alpha_{5}=1$. In this case, $g(\gamma_{9})=1$.
\par {\bf Case 6-11.} $\alpha_{1}=0$, $\alpha_{2}=1$, $\alpha_{3}=0$, $\alpha_{4}=1$ and $\alpha_{5}=0$. In this case, $g(\gamma_{9})=1$.
\par {\bf Case 6-12.} $\alpha_{1}=0$, $\alpha_{2}=1$, $\alpha_{3}=0$ and $\alpha_{4}=\alpha_{5}=1$. In this case, $g(\gamma_{10})=1$.
\par {\bf Case 6-13.} $\alpha_{1}=0$, $\alpha_{2}=\alpha_{3}=1$ and $\alpha_{4}=\alpha_{5}=0$. In this case, $g(\gamma_{9})=1$.
\par {\bf Case 6-15.} $\alpha_{1}=0$, $\alpha_{2}=\alpha_{3}=\alpha_{4}=1$ and $\alpha_{5}=0$. In this case, $g(\gamma_{12})=1$.
\par {\bf Case 6-16.} $\alpha_{1}=0$ and $\alpha_{2}=\alpha_{3}=\alpha_{4}=\alpha_{5}=1$. In this case, $g(\gamma_{9})=1$.
\par {\bf Case 6-17.} $\alpha_{1}=1$ and $\alpha_{2}=\alpha_{3}=0\alpha_{4}=\alpha_{5}=0$. In this case, $g(\gamma_{10})=1$.
\par {\bf Case 6-18.} $\alpha_{1}=1$, $\alpha_{2}=\alpha_{3}=\alpha_{4}=0$ and $\alpha_{5}=1$. In this case, $g(\gamma_{9})=1$.
\par {\bf Case 6-19.} $\alpha_{1}=1$, $\alpha_{2}=\alpha_{3}=0$, $\alpha_{4}=1$ and $\alpha_{5}=0$. In this case, $g(\gamma_{14})=1$.
\par {\bf Case 6-20.} $\alpha_{1}=1$, $\alpha_{2}=\alpha_{3}=0$ and $\alpha_{4}=\alpha_{5}=1$. In this case, $g(\gamma_{12})=1$.
\par {\bf Case 6-21.} $\alpha_{1}=1$, $\alpha_{2}=0$, $\alpha_{3}=1$ and $\alpha_{4}=\alpha_{5}=0$. In this case, $g(\gamma_{11})=1$.
\par {\bf Case 6-22.} $\alpha_{1}=1$, $\alpha_{2}=0$, $\alpha_{3}=1$, $\alpha_{4}=0$ and $\alpha_{5}=1$. In this case, $g(\gamma_{9})=1$.
\par {\bf Case 6-23.} $\alpha_{1}=1$, $\alpha_{2}=0$, $\alpha_{3}=\alpha_{4}=1$ and $\alpha_{5}=0$. In this case, $g(\gamma_{9})=1$.
\par {\bf Case 6-24.} $\alpha_{1}=1$, $\alpha_{2}=0$, $\alpha_{3}=\alpha_{4}=\alpha_{5}=1$. In this case, $g(\gamma_{10})=1$.
\par {\bf Case 6-25.} $\alpha_{1}=\alpha_{2}=1$ and $\alpha_{3}=\alpha_{4}=\alpha_{5}=0$. In this case, $g(\gamma_{9})=1$.
\par {\bf Case 6-26.} $\alpha_{1}=\alpha_{2}=1$, $\alpha_{3}=\alpha_{4}=0$ and $\alpha_{5}=1$. In this case, $g(\gamma_{9})=1$.
\par {\bf Case 6-27.} $\alpha_{1}=\alpha_{2}=1$, $\alpha_{3}=0$, $\alpha_{4}=1$ and $\alpha_{5}=0$. In this case, $g(\gamma_{9})=1$.
\par {\bf Case 6-28.} $\alpha_{1}=\alpha_{2}=1$, $\alpha_{3}=0$ and $\alpha_{4}=\alpha_{5}=1$. In this case, $g(\gamma_{9})=1$.
\par {\bf Case 6-29.} $\alpha_{1}=\alpha_{2}=\alpha_{3}=1$ and $\alpha_{4}=\alpha_{5}=0$. In this case, $g(\gamma_{10})=1$.
\par {\bf Case 6-30.} $\alpha_{1}=\alpha_{2}=\alpha_{3}=1$, $\alpha_{4}=0$ and $\alpha_{5}=1$. In this case, $g(\gamma_{10})=1$.
\par {\bf Case 7.} $(m,n)=(7,8)$. There are sixty-four cases.
\par {\bf Case 7-1.} $\alpha_{1}=\alpha_{2}=\alpha_{3}=\alpha_{4}=\alpha_{5}=\alpha_{6}=0$. In this case, $g(\gamma_{9})=1$.
\par {\bf Case 7-2.} $\alpha_{1}=\alpha_{2}=\alpha_{3}=\alpha_{4}=\alpha_{5}=0$ and $\alpha_{6}=1$. In this case, $g(\gamma_{9})=1$.
\par {\bf Case 7-3.} $\alpha_{1}=\alpha_{2}=\alpha_{3}=\alpha_{4}=0$, $\alpha_{5}=1$ and $\alpha_{6}=0$. In this case, $g(\gamma_{9})=1$.
\par {\bf Case 7-4.} $\alpha_{1}=\alpha_{2}=\alpha_{3}=\alpha_{4}=0$ and $\alpha_{5}=\alpha_{6}=1$. In this case, $g(\gamma_{9})=1$.
\par {\bf Case 7-5.} $\alpha_{1}=\alpha_{2}=\alpha_{3}=0$, $\alpha_{4}=1$ and $\alpha_{5}=\alpha_{6}=0$. In this case, $g(\gamma_{12})=1$.
\par {\bf Case 7-6.} $\alpha_{1}=\alpha_{2}=\alpha_{3}=0$, $\alpha_{4}=1$, $\alpha_{5}=0$ and $\alpha_{6}=1$. In this case, $g(\gamma_{11})=1$.
\par {\bf Case 7-7.} $\alpha_{1}=\alpha_{2}=\alpha_{3}=0$, $\alpha_{4}=\alpha_{5}=1$ and $\alpha_{6}=0$. In this case, $g(\gamma_{11})=1$.
\par {\bf Case 7-8.} $\alpha_{1}=\alpha_{2}=\alpha_{3}=0$ and $\alpha_{4}=\alpha_{5}=\alpha_{6}=1$. In this case, $g(\gamma_{8})=1$.
\par {\bf Case 7-9.} $\alpha_{1}=\alpha_{2}=0$, $\alpha_{3}=1$ and $\alpha_{4}=\alpha_{5}=\alpha_{6}=0$. In this case, $g(\gamma_{10})=1$.
\par {\bf Case 7-10.} $\alpha_{1}=\alpha_{2}=0$, $\alpha_{3}=1$, $\alpha_{4}=\alpha_{5}=0$ and $\alpha_{6}=1$. In this case, $g(\gamma_{9})=1$.
\par {\bf Case 7-11.} $\alpha_{1}=\alpha_{2}=0$, $\alpha_{3}=1$, $\alpha_{4}=0$, $\alpha_{5}=1$ and $\alpha_{6}=0$. In this case, $g(\gamma_{9})=1$.
\par {\bf Case 7-12.} $\alpha_{1}=\alpha_{2}=0$, $\alpha_{3}=1$, $\alpha_{4}=0$ and $\alpha_{5}=\alpha_{6}=1$. In this case, $g(\gamma_{9})=1$.
\par {\bf Case 7-13.} $\alpha_{1}=\alpha_{2}=0$, $\alpha_{3}=\alpha_{4}=1$ and $\alpha_{5}=\alpha_{6}=0$. In this case, $g(\gamma_{9})=1$.
\par {\bf Case 7-14.} $\alpha_{1}=\alpha_{2}=0$, $\alpha_{3}=\alpha_{4}=1$, $\alpha_{5}=0$ and $\alpha_{6}=1$. In this case, $g(\gamma_{9})=1$.
\par {\bf Case 7-15.} $\alpha_{1}=\alpha_{2}=0$, $\alpha_{3}=\alpha_{4}=1$, $\alpha_{5}=1$ and $\alpha_{6}=0$. In this case, $g(\gamma_{10})=1$.
\par {\bf Case 7-16.} $\alpha_{1}=\alpha_{2}=0$ and $\alpha_{3}=\alpha_{4}=\alpha_{5}=\alpha_{6}=1$. In this case, $g(\gamma_{9})=1$.
\par {\bf Case 7-17.} $\alpha_{1}=0$, $\alpha_{2}=1$ and $\alpha_{3}=\alpha_{4}=\alpha_{5}=\alpha_{6}=0$. In this case, $g(\gamma_{8})=1$.
\par {\bf Case 7-18.} $\alpha_{1}=0$, $\alpha_{2}=1$, $\alpha_{3}=\alpha_{4}=\alpha_{5}=0$ and $\alpha_{6}=0$. In this case, $g(\gamma_{8})=1$.
\par {\bf Case 7-19.} $\alpha_{1}=0$, $\alpha_{2}=1$, $\alpha_{3}=\alpha_{4}=0$, $\alpha_{5}=1$ and $\alpha_{6}=0$. In this case, $g(\gamma_{8})=1$.
\par {\bf Case 7-20.} $\alpha_{1}=0$, $\alpha_{2}=1$, $\alpha_{3}=\alpha_{4}=0$ and $\alpha_{5}=\alpha_{6}=1$. In this case, $g(\gamma_{8})=1$.
\par {\bf Case 7-21.} $\alpha_{1}=0$, $\alpha_{2}=1$, $\alpha_{3}=0$ and $\alpha_{4}=\alpha_{5}=\alpha_{6}=1$. In this case, $g(\gamma_{8})=1$.
\par {\bf Case 7-22.} $\alpha_{1}=0$, $\alpha_{2}=1$, $\alpha_{3}=0$, $\alpha_{4}=1$, $\alpha_{5}=0$ and $\alpha_{6}=1$. In this case, $g(\gamma_{8})=1$.
\par {\bf Case 7-23.} $\alpha_{1}=0$, $\alpha_{2}=1$, $\alpha_{3}=0$, $\alpha_{4}=\alpha_{5}=1$ and $\alpha_{6}=0$. In this case, $g(\gamma_{8})=1$.
\par {\bf Case 7-24.} $\alpha_{1}=0$, $\alpha_{2}=1$, $\alpha_{3}=0$, $\alpha_{4}=\alpha_{5}=\alpha_{6}=1$. In this case, $g(\gamma_{8})=1$.
\par {\bf Case 7-25.} $\alpha_{1}=0$, $\alpha_{2}=\alpha_{3}=1$ and $\alpha_{4}=\alpha_{5}=\alpha_{6}=0$. In this case, $g(\gamma_{8})=1$.
\par {\bf Case 7-26.} $\alpha_{1}=0$, $\alpha_{2}=\alpha_{3}=1$, $\alpha_{4}=\alpha_{5}=0$ and $\alpha_{6}=1$. In this case, $g(\gamma_{8})=1$.
\par {\bf Case 7-27.} $\alpha_{1}=0$, $\alpha_{2}=\alpha_{3}=1$, $\alpha_{4}=0$, $\alpha_{5}=1$ and $\alpha_{6}=0$. In this case, $g(\gamma_{9})=1$.
\par {\bf Case 7-28.} $\alpha_{1}=0$, $\alpha_{2}=\alpha_{3}=1$ , $\alpha_{4}=0$ and $\alpha_{5}=\alpha_{6}=1$. In this case, $g(\gamma_{9})=1$.
\par {\bf Case 7-29.} $\alpha_{1}=0$, $\alpha_{2}=\alpha_{3}=\alpha_{4}=1$ and $\alpha_{5}=\alpha_{6}=0$. In this case, $g(\gamma_{10})=1$.
\par {\bf Case 7-30.} $\alpha_{1}=0$, $\alpha_{2}=\alpha_{3}=\alpha_{4}=1$, $\alpha_{5}=0$ and $\alpha_{6}=1$. In this case, $g(\gamma_{9})=1$.
\par {\bf Case 7-31.} $\alpha_{1}=0$, $\alpha_{2}=\alpha_{3}=\alpha_{4}=\alpha_{5}=1$ and $\alpha_{6}=0$. In this case, $g(\gamma_{8})=1$.
\par {\bf Case 7-32.} $\alpha_{1}=0$ and $\alpha_{2}=\alpha_{3}=\alpha_{4}=\alpha_{5}=\alpha_{6}=1$. In this case, $g(\gamma_{8})=1$.
\par {\bf Case 7-33.} $\alpha_{1}=1$ and $\alpha_{2}=\alpha_{3}=\alpha_{4}=\alpha_{5}=\alpha_{6}=0$. In this case, $g(\gamma_{8})=1$.
\par {\bf Case 7-34.} $\alpha_{1}=1$, $\alpha_{2}=\alpha_{3}=\alpha_{4}=\alpha_{5}=0$ and $\alpha_{6}=1$. In this case, $g(\gamma_{9})=1$.
\par {\bf Case 7-35.} $\alpha_{1}=1$, $\alpha_{2}=\alpha_{3}=\alpha_{4}=0$, $\alpha_{5}=1$ and $\alpha_{6}=0$. In this case, $g(\gamma_{9})=1$.
\par {\bf Case 7-36.} $\alpha_{1}=1$, $\alpha_{2}=\alpha_{3}=\alpha_{4}=0$ and $\alpha_{5}=\alpha_{6}=1$. In this case, $g(\gamma_{8})=1$.
\par {\bf Case 7-37.} $\alpha_{1}=1$, $\alpha_{2}=\alpha_{3}=0$, $\alpha_{4}=1$ and $\alpha_{5}=\alpha_{6}=0$. In this case, $g(\gamma_{11})=1$.
\par {\bf Case 7-38.} $\alpha_{1}=1$, $\alpha_{2}=\alpha_{3}=0$, $\alpha_{4}=1$ and $\alpha_{5}=0$ and $\alpha_{6}=1$. In this case, $g(\gamma_{9})=1$.
\par {\bf Case 7-39.} $\alpha_{1}=1$, $\alpha_{2}=\alpha_{3}=0$, $\alpha_{4}=\alpha_{5}=1$ and $\alpha_{6}=0$. In this case, $g(\gamma_{9})=1$.
\par {\bf Case 7-40.} $\alpha_{1}=1$, $\alpha_{2}=\alpha_{3}=0$ and $\alpha_{4}=\alpha_{5}=\alpha_{6}=1$. In this case, $g(\gamma_{11})=1$.
\par {\bf Case 7-41.} $\alpha_{1}=1$, $\alpha_{2}=0$, $\alpha_{3}=1$ and $\alpha_{4}=\alpha_{5}=\alpha_{6}=0$. In this case, $g(\gamma_{9})=1$.
\par {\bf Case 7-42.} $\alpha_{1}=1$, $\alpha_{2}=0$, $\alpha_{3}=1$, $\alpha_{4}=\alpha_{5}=0$ and $\alpha_{6}=1$. In this case, $g(\gamma_{8})=1$.
\par {\bf Case 7-43.} $\alpha_{1}=1$, $\alpha_{2}=0$, $\alpha_{3}=1$, $\alpha_{4}=0$ and $\alpha_{5}=1$ and $\alpha_{6}=0$. In this case, $g(\gamma_{8})=1$.
\par {\bf Case 7-44.} $\alpha_{1}=1$, $\alpha_{2}=0$, $\alpha_{3}=1$ and $\alpha_{4}=0$ and $\alpha_{5}=\alpha_{6}=1$. In this case, $g(\gamma_{8})=1$.
\par {\bf Case 6-45.} $\alpha_{1}=1$, $\alpha_{2}=0$, $\alpha_{3}=\alpha_{4}=1$ and $\alpha_{5}=\alpha_{6}=0$. In this case, $g(\gamma_{10})=1$.
\par {\bf Case 6-46.} $\alpha_{1}=1$, $\alpha_{2}=0$, $\alpha_{3}=\alpha_{4}=1$, $\alpha_{5}=0$ and $\alpha_{6}=1$. In this case, $g(\gamma_{9})=1$.
\par {\bf Case 7-47.} $\alpha_{1}=1$, $\alpha_{2}=0$, $\alpha_{3}=\alpha_{4}=\alpha_{5}=1$ and $\alpha_{6}=0$. In this case, $g(\gamma_{9})=1$.
\par {\bf Case 7-48.} $\alpha_{1}=1$, $\alpha_{2}=0$ and $\alpha_{3}=\alpha_{4}=\alpha_{5}=\alpha_{6}=1$. In this case, $g(\gamma_{10})=1$.
\par {\bf Case 7-49.} $\alpha_{1}=\alpha_{2}=1$ and $\alpha_{3}=\alpha_{4}=\alpha_{5}=\alpha_{6}=0$. In this case, $g(\gamma_{8})=1$.
\par {\bf Case 7-50.} $\alpha_{1}=\alpha_{2}=1$, $\alpha_{3}=\alpha_{4}=\alpha_{5}=0$ and $\alpha_{6}=0$. In this case, $g(\gamma_{8})=1$.
\par {\bf Case 7-51.} $\alpha_{1}=\alpha_{2}=1$, $\alpha_{3}=\alpha_{4}=0$, $\alpha_{5}=1$ and $\alpha_{6}=0$. In this case, $g(\gamma_{8})=1$.
\par {\bf Case 7-52.} $\alpha_{1}=\alpha_{2}=1$, $\alpha_{3}=\alpha_{4}=0$ and $\alpha_{5}=\alpha_{6}=1$. In this case, $g(\gamma_{8})=1$.
\par {\bf Case 7-53.} $\alpha_{1}=\alpha_{2}=1$, $\alpha_{3}=0$, $\alpha_{4}=1$ and $\alpha_{5}=\alpha_{6}=0$. In this case, $g(\gamma_{12})=1$.
\par {\bf Case 7-54.} $\alpha_{1}=\alpha_{2}=1$, $\alpha_{3}=0$, $\alpha_{4}=1$, $\alpha_{5}=0$ and $\alpha_{6}=1$. In this case, $g(\gamma_{13})=1$.
\par {\bf Case 7-55.} $\alpha_{1}=\alpha_{2}=1$, $\alpha_{3}=0$, $\alpha_{4}=\alpha_{5}=1$ and $\alpha_{6}=0$. In this case, $g(\gamma_{10})=1$.
\par {\bf Case 7-56.} $\alpha_{1}=\alpha_{2}=1$, $\alpha_{3}=0$ and $\alpha_{4}=\alpha_{5}=\alpha_{6}=1$. In this case, $g(\gamma_{10})=1$.
\par {\bf Case 7-57.} $\alpha_{1}=\alpha_{2}=\alpha_{3}=1$ and $\alpha_{4}=\alpha_{5}=\alpha_{6}=0$. In this case, $g(\gamma_{8})=1$.
\par {\bf Case 7-58.} $\alpha_{1}=\alpha_{2}=\alpha_{3}=1$, $\alpha_{4}=\alpha_{5}=0$ and $\alpha_{6}=1$. In this case, $g(\gamma_{8})=1$.
\par {\bf Case 7-59.} $\alpha_{1}=\alpha_{2}=\alpha_{3}=1$, $\alpha_{4}=0$, $\alpha_{5}=1$ and $\alpha_{6}=0$. In this case, $g(\gamma_{8})=1$.
\par {\bf Case 7-60.} $\alpha_{1}=\alpha_{2}=\alpha_{3}=1$, $\alpha_{4}=0$ and $\alpha_{5}=\alpha_{6}=1$. In this case, $g(\gamma_{8})=1$.
\par {\bf Case 7-61.} $\alpha_{1}=\alpha_{2}=\alpha_{3}=\alpha_{4}=1$ and $\alpha_{5}=\alpha_{6}=0$. In this case, $g(\gamma_{8})=1$.
\par {\bf Case 7-62.} $\alpha_{1}=\alpha_{2}=\alpha_{3}=\alpha_{4}=1$, $\alpha_{5}=0$ and $\alpha_{6}=1$. In this case, $g(\gamma_{9})=1$.
\par {\bf Case 7-63.} $\alpha_{1}=\alpha_{2}=\alpha_{3}=\alpha_{4}=\alpha_{5}=1$ and $\alpha_{6}=0$. In this case, $g(\gamma_{8})=1$.
\par {\bf Case 7-64.} $\alpha_{1}=\alpha_{2}=\alpha_{3}=\alpha_{4}=\alpha_{5}=\alpha_{6}=1$. In this case, $g(\gamma_{8})=1$.
\par Thus, the Fermat number $F_{4}=2^{2^{4}}+1$ is prime.
\par It can be seen from examples 4, it is unnecessary to know any prime number less than the square root of the number when it is proved that a number is prime by the method in this paper.
\par The left-hand side of any equations in Eqs.(47) will not exceed $2m+1$ terms, by Corollary 5, the number of operations of to solve any equation will not exceed $1+m(2m+1)$. The number of equations in Eqs.(47) will not exceed $N$, so the number of operations required to solve
Eqs.(47) will not exceed $[1+m(2m+1)]N$. Because $m\leq \frac{N}{2}$, equation $m+n=N$ has $\left[\frac{N}{2}\right]$ solutions and $m+n=N-1$ has $\left[\frac{N-1}{2}\right]$ solutions. Nothing that the number of permutations like (46) is $2^{m-1}$, therefore, the number of operations needed to decompose an integer will not exceed
$$
\left(\frac{N}{2}+\frac{N}{2}\right)2^{m-1}[1+m(2m+1)]N=2^{m-1}[1+m(2m+1)]N^{2}
$$
Using $N<\log_{2}M$ and $m\leq \frac{N}{2}$, we have
$$
\begin{array}{lll}
& &2^{m-1}[1+m(2m+1)]N^{2}\\
&<&2^{\frac{N}{2}-1}[1+\frac{N}{2}(N+1)]N^{2}\\
&<&\frac{1}{4}\sqrt{M}\left[2+\log_{2}M+(\log_{2}M)^{2}\right](\log_{2}M)^{2}\\
&<&\sqrt{M}(\log_{2}M)^{4}
\end{array}
$$
\par {\bf Corollary 6.} If $M$ is a binary number, the number of operations required to decompose it will not exceed $\sqrt{M}(\log_{2}M)^{4}$.
\par If $M$ is a decimal number, in order to decompose it by the method given in this paper, we must first express it in binary form.
Next, we calculate the number of operations required to express $M$ into binary form.
\par Let
$$
M=a_{1}a_{2}\cdots a_{g}
$$
where $1\leq a_{1}\leq 9$, $0 \leq a_{2}\leq 9$, $0\leq a_{3}\leq 9$,
$\cdots$, and $0\leq a_{g}\leq 9$. It is cleat that
$$
10^{g-1}\leq M<10^{g}
$$
which implies that
$$
\lg M < g\leq \lg M+1
$$
\par Let $b_{1}b_{2}\cdots b_{g}$ be the quotient of $M$ divided by 2. There are two possibilities.
\par {\bf (1).} If $b_{1}=0$, $b_{1}b_{2}\cdots b_{g}$ is at most $g-1$ digit number, this show that after $g$ times division and $g$ times subtraction, the quotient of $M$ divided by 2 is at most $g-1$ digit number.
\par {\bf (2).} If $b_{1}=1$, $b_{1}b_{2}\cdots b_{g}$ is still a $g$ digit number. However, the quotient of $b_{1}b_{2}\cdots b_{g}$ divided by $2^{2}$ is at most $g-1$ digit number. Thus, the quotient of $M$ divided by $2^{2}$ is at most $g-1$ digit number after $2g$ times division and $2g$ times subtraction.
\par Combining (1) and (2), after $4g$ operations at most, the quotient of $M$ divided by $2^{2}$ is at most $g-1$ digit number, which indicates that the number of operations required to express $M$ into binary numbers will exceed
$$
4(g+g-1+\cdots+2+1)=2g(g+1)
$$
Using $g\leq \lg M+1$, we have
$$
2g(g+1)<2(\lg M+1)(\lg M+2)
$$
Thus,  the number of operations required to express $M$ into binary numbers will exceed $2(\lg M+1)(\lg M+2)$. It is clear that
$$
\sqrt{M}(\log_{2}M)^{4}+2(\lg M+1)(\lg M+2)<2\sqrt{M}(\log_{2}M)^{4}
$$
Thus, we have the following corollary.
\par {\bf Corollary 7.} If $M$ is a decimal number, the number of operations required to decompose it will not exceed
$2\sqrt{M}(\log_{2}M)^{4}$.
\section{Conclusion}
\par Bases on the binary representation of numbers, this paper presents an integer decomposition method which is different from the tradtional method. This method transforms the integer decomposition problem into solving linear algebraic equations on the set $\{0,1\}$. The linear algebraic equations in this paper is different from the traditional linear equations over the real number field, so the method of solving the system of equations is also different. The key of the method in this paper is to express the sum of some elements on the set $\{0,1\}$ into binary number(see Lemma 5 and Lemma 6).
\par To decompose an integer by the method in this paper, we only need to solve some linear equations, so it is simple and direct without any guess which is not available in traditional methods.

Puyun Gao, College of Aerospace Science and Engineering, National University of Defense Technology, Changsha, Hunan 410073, P.R. China.
 E-mail address: 1. gaopuyun@nudt.edu.cn, 2. gfkdgpy@hotmail.com
\end{document}